\documentclass[reqno,11pt]{amsart}
\usepackage[foot]{amsaddr}
\usepackage{amssymb,amsmath,amsthm,amsfonts}
\usepackage{mathrsfs,dsfont,comment,mathscinet,mathtools}
\usepackage{regexpatch}
\usepackage{graphicx,pgfplots}
\usepackage{bigstrut}

\usepackage[utf8]{inputenc}
\usepackage[normalem]{ulem}
\usepackage[left=2.5cm,right=2.5cm,top=2.5cm,bottom=2.5cm]{geometry}

\usepackage{pifont}
\usepackage{tkz-fct}

\usepackage{graphicx,tikz,color}
\usetikzlibrary{quotes,angles}

\usepackage[format=hang,labelfont=bf]{caption}

\usepackage{enumitem}

\usepackage{esint}

\usepackage{ wasysym }

\usepackage{epstopdf}
\usepackage[colorlinks=true, pdfstartview=FitV, linkcolor=blue, 
            citecolor=blue, urlcolor=blue]{hyperref}
            
\allowdisplaybreaks
\numberwithin{equation}{section}
\theoremstyle{plain}
\newtheorem{theorem}{Theorem}[section]
\newtheorem{proposition}[theorem]{Proposition}
\newtheorem{lemma}[theorem]{Lemma}

\newtheorem{corollary}[theorem]{Corollary}
\newtheorem{definition}[theorem]{Definition}
\theoremstyle{definition}

\newtheorem*{theorem*}{Theorem}

\definecolor{mblue}{HTML}{13439b}

\newcommand{\N}{\mathbb{N}}
\newcommand{\R}{\mathbb{R}}

\newcommand{\RNN}{\mathbb{R}^{N\times N}}

\newcommand{\Rt}{\mathbb{R}^3}
\newcommand{\Rtt}{\mathbb{R}^{3\times 3}}
\newcommand{\St}{\mathbb{S}^2}

\renewcommand{\d}{\mathrm{d}}
\newcommand{\restr}[1]{|_{#1}}

\def\Xint#1{\mathchoice
{\XXint\displaystyle\textstyle{#1}}%
{\XXint\textstyle\scriptstyle{#1}}%
{\XXint\scriptstyle\scriptscriptstyle{#1}}%
{\XXint\scriptscriptstyle\scriptscriptstyle{#1}}%
\!\int}
\def\XXint#1#2#3{{\setbox0=\hbox{$#1{#2#3}{\int}$ }
\vcenter{\hbox{$#2#3$ }}\kern-.6\wd0}}

\def\dashint{\Xint-}

\makeatletter
\newcommand{\subalign}[1]{%
	\vcenter{%
		\Let@ \restore@math@cr \default@tag
		\baselineskip\fontdimen10 \scriptfont\tw@
		\advance\baselineskip\fontdimen12 \scriptfont\tw@
		\lineskip\thr@@\fontdimen8 \scriptfont\thr@@
		\lineskiplimit\lineskip
		\ialign{\hfil$\m@th\scriptstyle##$&$\m@th\scriptstyle{}##$\hfil\crcr
			#1\crcr
		}%
	}%
}
\makeatother

\makeatletter
\newcommand{\subsetcong}{\mathrel{\mathpalette\subseteq@cong\relax}}
\newcommand{\subsetsim}{\mathrel{\mathpalette\subset@sim\relax}}

\newcommand{\subseteq@cong}[2]{%
	\vbox{\offinterlineskip\m@th
		\ialign{\hfil$#1##$\hfil\cr
			\sim\cr\subseteq\cr
		}%
	}%
}

\newcommand{\subset@sim}[2]{%
	\vbox{\offinterlineskip\m@th
		\ialign{\hfil$#1##$\hfil\cr
			\sim\cr\subset\cr
		}%
	}%
}

\makeatother

\newcommand{\adj}{\mathrm{adj}}
\newcommand{\sym}{\mathrm{sym}}

\newcommand{\per}{\mathrm{Per}}
\renewcommand{\div}{\mathrm{div}}
\newcommand{\curl}{\mathrm{curl}}
\newcommand{\dist}{\mathrm{dist}}
\newcommand\wk{\rightharpoonup}

\newcommand*\closure[1]{\overline{#1}}

\newcommand{\leb}{\mathscr{L}^3}
\newcommand{\haus}{\mathscr{H}^{2}}

\newcommand{\imt}{\mathrm{im}_{\rm T}}
\newcommand{\img}{\mathrm{im}_{\rm G}}
\newcommand{\domg}{\mathrm{dom}_{\rm G}}
\newcommand{\Det}{\mathrm{Det}}

\newcommand\EEE{\color{black}}

\newcommand{\MMM}{\color{black}} 

\makeatletter
\xpatchcmd{\proof}{\itshape}{\bfseries}{}{}
\makeatother

\title[Domain structures in magnetostrictive solids]{A variational approach to the emergence of domain structures in magnetostrictive solids}
\author[M. Bresciani]{Marco Bresciani${}^{*}$}
\address{* Department of Mathematics, Friedrich-Alexander-Universit\"{a}t Erlangen-N\"{u}rnberg, Cauerstrasse 11, 91058 Erlangen, Germany}
\email{marco.bresciani@fau.de}
\author[M. Friedrich]{Manuel Friedrich${}^{*}$}
\email{manuel.friedrich@fau.de}

\date{\today}

\begin{document}

\setlength\parindent{0pt}

\vskip .2truecm
\begin{abstract}
We consider a variational model for magnetoelastic solids in the large-strain setting with the magnetization field defined on the unknown deformed configuration. Through a simultaneous linearization of the deformation and sharp-interface limit of the magnetization with respect to the easy axes, performed in terms of $\Gamma$-convergence, we identify a variational model describing the formation of domain structures and accounting for magnetostriction in the small-strain setting. Our analysis incorporates the effect of  magnetic fields and, for uniaxial materials, ensures the regularity of minimizers of the effective energy. 
\end{abstract}
\maketitle

\section{Introduction}
\label{sec:intro}

\subsection{Motivation}
In the modeling of deformable ferromagnets, magnetostriction refers to the effect for  which mechanical and magnetic response of materials influence each other. Nowadays, this effect is  largely exploited in the design of many technological devices such as sensors and actuators (see, e.g., \cite{actuators,sensors}). In this regard, specific rare-earth and shape-memory alloys (e.g., TerFeNOL and NiMnGa)   have been found of great  interest in view of their ability of undergoing remarkably large deformations upon the application of relatively moderate magnetic field, which is often referred to as giant magnetostriction \cite{engdhal}.

Magnetocrystalline anisotropy, i.e., the existence of preferred  directions  of magnetization  (called easy axes)  determined by the underlying lattice structure of materials, plays a key role in the mechanism behind  magnetostriction. In absence of external stimuli, magnetic materials spontaneously  arrange  themselves in regions of uniform magnetization  (called Weiss domains) \EEE separated by thin transition layers (called Bloch walls) \EEE \cite{hubert.schaefer}. The application of an external magnetic field leads to a reorganization of the domain structure mainly due to the more favorable orientation of  certain  easy  axes \EEE with respect to the  magnetic field.  The movement of domains results in a macroscopic strain, thus inducing a deformation of the material. For a comprehensive treatment of the physics of magnetoelastic interactions, we refer, e.g., to \cite{brown,landau.lifshitz}.

Providing an effective description for the emergence of  domain structures accounting for mechanical stresses in precise mathematical terms appears to be a very challenging task. This study  constitutes
 \EEE an attempt in this direction by  means of variational methods. Starting from a large-strain model for magnetostrictive solids, we study the asymptotic behavior of the energy for infinitesimal strains and for magnetization aligning with the easy axes. In this way, we  identify an effective model accounting for the linearized elastic energy and the formation of domain  structures. \EEE Our main results are outlined in Subsection \ref{subsec:main} below.

 \EEE 

\subsection{The variational model} \label{subsec:brown}
In this work, we adopt the  variational model of magnetoelasticity due to Brown \cite{brown,desimone.podioguidugli}.
A magnetoelastic body is subjected to a deformation $\boldsymbol{y}\colon \Omega \to \Rt$, where $\Omega\subset \Rt$ represents its reference configuration, and  a magnetization field $\boldsymbol{m}\colon \boldsymbol{y}(\Omega)\to \Rt$. The latter is understood as the average of magnetic dipoles per unit volume. Hence, the saturation constraint (see, e.g., \cite[p.~73]{brown}) requires  $|\boldsymbol{m}|=1$ in $\boldsymbol{y}(\Omega)$.
The magnetoelastic energy reads
\begin{equation}
	\label{eq:intro:B}
	\begin{split}
		(\boldsymbol{y},\boldsymbol{m})\mapsto & \int_\Omega \Psi(D\boldsymbol{y},\boldsymbol{m}\circ \boldsymbol{y})\,\d\boldsymbol{x}
		+\int_{\boldsymbol{y}(\Omega)}  \Phi \left( \left( (D\boldsymbol{y})^\top \circ \boldsymbol{y}^{-1}\right ) \boldsymbol{m}   \right)  \,\d\boldsymbol{\xi}\\
		&+ \int_{\boldsymbol{y}(\Omega)} |D\boldsymbol{m}|^2\,\d\boldsymbol{\xi}+ \int_{\Rt} |\boldsymbol{h}^{\boldsymbol{y},\boldsymbol{m}}|^2\,\d\boldsymbol{\xi}.
	\end{split}
\end{equation}
The first term in \eqref{eq:intro:B} accounts for the elastic energy and involves a material density $\Psi$ which exhibits the coupling   between magnetization and deformation gradient. The prototypical example (see, e.g., \cite{james.kinderlehrer}) takes the form
\begin{equation}
	\label{eq:Psi}
	\Psi(\boldsymbol{F},\boldsymbol{z})=W(\boldsymbol{\Lambda}^{-1}(\boldsymbol{z})\boldsymbol{F}),
\end{equation} 
where $W$ is an auxiliary density and $\boldsymbol{\Lambda}(\boldsymbol{z})$  represents  the spontaneous strain corresponding to $\boldsymbol{z}$.
The second term in \eqref{eq:intro:B}  stands  for the anisotropy energy. The density $\Phi$ is nonnegative and vanishes exactly in the direction of the easy axes, thus favoring their attainment. For instance, in the case of uniaxial materials, 
\begin{equation}
	\label{eq:uniaxial}
	  \Phi(\boldsymbol{z})=\kappa \left( 1-(\boldsymbol{a}\cdot \boldsymbol{z})^2 \right)
\end{equation}
with  $\boldsymbol{a}\in\St$ being the easy axis and $\kappa>0$, while  for cubic materials, 
\begin{equation}
	\label{eq:cubic}
	\begin{split}
		\Phi(\boldsymbol{z})=&   \kappa_1 \left( (\boldsymbol{z}\cdot \boldsymbol{a}_1)^2(\boldsymbol{z}\cdot \boldsymbol{a}_2)^2 + (\boldsymbol{z}\cdot \boldsymbol{a}_1)^2(\boldsymbol{z}\cdot \boldsymbol{a}_3)^2 + (\boldsymbol{z}\cdot \boldsymbol{a}_2)^2(\boldsymbol{z}\cdot \boldsymbol{a}_3)^2   \right)\\
		&+ \kappa_2 (\boldsymbol{z}\cdot \boldsymbol{a}_1)^2(\boldsymbol{z}\cdot \boldsymbol{a}_2)^2(\boldsymbol{z}\cdot \boldsymbol{a}_3)^2,
	\end{split}
\end{equation}
where $\boldsymbol{a}_1,\boldsymbol{a}_2,\boldsymbol{a}_3 \in\St$ are the   mutually orthogonal   easy axes, $\kappa_1>0$, and $\kappa_2\geq 0$.  
The third term in \eqref{eq:intro:B} constitutes the exchange energy which  favors  a constant magnetization of the specimen. More generally, the norm squared  could  be replaced by a  quadratic form in the magnetization gradient. Eventually, the fourth term in \eqref{eq:intro:B} stands for the magnetostatic self-energy related to the long-range interactions. In absence of external currents, the stray field $\boldsymbol{h}^{\boldsymbol{y},\boldsymbol{m}}$ solves the Maxwell system
\begin{equation*}
	\begin{cases}
		\mathrm{curl}\, \boldsymbol{h}^{\boldsymbol{y},\boldsymbol{m}}=\boldsymbol{0} & {}\\
		\mathrm{div}\,\boldsymbol{h}^{\boldsymbol{y},\boldsymbol{m}}=-\mathrm{div} \left( \chi_{\boldsymbol{y}(\Omega)} \boldsymbol{m}\det D \boldsymbol{y}^{-1} \right) & {}
	\end{cases}
	\quad \text{in $\Rt$.}
\end{equation*}
We remark that, when $\Psi$ is invariant under left-multiplication by rotation matrices (see, e.g., \cite[Eq. (20)]{james.kinderlehrer}),  the energy in \eqref{eq:intro:B} fulfills the principle of frame indifference. In particular, this happens for \eqref{eq:Psi} when the same property is satisfied by $W$  and $\boldsymbol{\Lambda}$ is an objective tensor field. 
 In this regard, we highlight  that the  dependence of the anisotropy energy on the deformation gradient   (see, e.g., \cite[Eq.~(21)]{james.kinderlehrer}) remains crucial for \EEE  frame indifference. Additionally, the energy is invariant upon replacing $\boldsymbol{m}$ with $-\boldsymbol{m}$.  

Equilibrium configurations correspond to minimizers of the energy \eqref{eq:intro:B}. The minimization becomes nontrivial once the deformation $\boldsymbol{y}$ is prescribed on a portion of the boundary or applied loads are considered. In presence of an external field $\boldsymbol{f}$, the corresponding energy contribution reads
\begin{equation*}
	(\boldsymbol{y},\boldsymbol{m})\mapsto \int_{\boldsymbol{y}(\Omega)} \boldsymbol{f}\cdot \boldsymbol{m}\,\d\boldsymbol{\xi}.
\end{equation*}

\subsection{Review of the most relevant literature}

To put our findings into context, we briefly review the most relevant literature. Without any claim of completeness, we only mention a few works concerning the model outlined in Subsection \ref{subsec:brown} that are related to our work.

The existence of minimizers for the functional \eqref{eq:intro:B} has been investigated in  \cite{barchiesi.henao.moracorral,bresciani,bresciani.davoli.kruzik, bresciani.friedrich.moracorral, kruzik.stefanelli.zeman}   under  various \EEE regularity  assumptions on admissible deformations.    The works \cite{bresciani,bresciani.davoli.kruzik,kruzik.stefanelli.zeman} address also the   existence of quasistatic evolution in the rate-independent \EEE setting. Other results of this kind  have been achieved in \cite{barchiesi.desimone,bresciani.stroffolini,henao.stroffolini} for a similar model of nematic elastomers. Note that in all these works the anisotropy energy is either  assumed to be independent on the deformation gradient \cite[Eq.~(1.2)]{bresciani.davoli.kruzik} or entirely neglected. Indeed, according to \eqref{eq:intro:B}, the anisotropy contribution constitutes a leading-order term of the functional whose lower semicontinuity cannot be guaranteed when $\Phi$ is nonconvex such as in \eqref{eq:uniaxial}--\eqref{eq:cubic}. 
 \EEE 

The linearization of the energy \eqref{eq:intro:B}  has been rigorously performed in \cite{almi.kruzik.molchanova},  thus recovering a well-established magnetoelasticity model for the small-strain setting \cite{james.desimone}.  For rigid bodies, the sharp-interface limit of magnetizations towards the easy axes has been addressed in \cite{anzellotti.baldo.visintin}.  The case of \MMM deformable \EEE bodies has been tackled in \cite{gkms-arma,gkms}. More precisely, these works concern   
 multiphase solids with Eulerian interfaces but,  as mentioned in \cite[p.~411]{gkms},  their analysis  can be reformulated in terms of the model in Subsection~\ref{subsec:brown}  by interpreting the diffuse phase indicators as magnetization fields. Still,  the dependence of the anisotropy energy on the deformation gradient is dropped in these two works, so that the  model does not comply with frame indifference. All above-mentioned results for the linearization and the sharp-interface limit are stated in terms of $\Gamma$-convergence \cite{braides,dalmaso}.   \EEE 

 Although not directly relevant for our study, we  also mention the derivation of effective models, still by $\Gamma$-convergence, \EEE in the direction  of \EEE dimension reduction \cite{bresciani.m3as,bk} and relaxation \cite{moracorral.oliva,scilla.stroffolini}  .

\subsection{Contributions of the paper} \label{subsec:main}
In this paper, we rigorously derive a variational model describing the formation of  domain \EEE structures and accounting for magnetostriction in the small-strain regime. The derivation is performed in the sense of $\Gamma$-convergence through a simultaneous linearization of deformations and sharp-interface limit of magnetizations with respect to the easy axes.
Our main results are given in Theorem \ref{thm:compactness}, Theorem \ref{thm:gamma}, and Theorem \ref{thm:regularity} below. Here, we limit ourselves in briefly  outlining \EEE their content and we refer to Section \ref{sec:main} for their precise statement.  

We adopt the  functional setting of \cite{barchiesi.henao.moracorral,bresciani}. Deformations are  modeled as maps $\boldsymbol{y}\in W^{1,p}(\Omega;\Rt)$ with $p>2$  \MMM satisfying the divergence identities (see Definition \ref{def:div} below), i.e., they \EEE do not exhibit cavitation. By replacing the deformed configuration $\boldsymbol{y}(\Omega)$ with the topological image (see Definition \ref{def:topim} below), we can define magnetizations as Sobolev maps on this set, which in turn allows for a meaningful definition of the magnetoelastic energy \eqref{eq:intro:B}.

For a small parameter $\varepsilon>0$, we consider a rescaled version $\mathcal{E}_\varepsilon$ of the magnetoelastic energy (see  \eqref{eq:Eeps}--\eqref{eq:Em-eps} below). Specifically, the elastic term is rescaled by $\varepsilon^2$ and the elastic density takes the form in \eqref{eq:Psi} with the spontaneous strain of magnetizations (i.e., the term $\boldsymbol{\Lambda}_\varepsilon(\boldsymbol{m}\circ \boldsymbol{y})$ in \eqref{eq:Ee-eps} below) given by a perturbation of the identity of order $\varepsilon$  (see \eqref{eq:Lambda} below); \EEE anisotropy and exchange term are rescaled similarly to the classical Modica-Mortola functional \cite{baldo}  in terms of a power $\varepsilon^\beta$   with
\begin{equation}
	\label{eq:beta}
	0<\beta < \frac{2(q-1)}{q}\wedge 1,
\end{equation}
where $q>1$ specifies the growth of $W$ with respect to  the determinant of its argument (see assumption (W5) below); \EEE  eventually,  we neglect  the stray-field term. \EEE Condition \eqref{eq:beta} is commented in Subsection \ref{subsec:discussion} below.

Theorem \ref{thm:compactness} shows that the sequence $(\mathcal{E}_\varepsilon)_{\varepsilon>0}$ enjoys suitable compactness properties or, in the terminology of $\Gamma$-convergence, that this sequence is equi-coercive (see, e.g., \cite[Definition 7.6]{dalmaso}). Limiting states are determined by infinitesimal displacements $\boldsymbol{u}\colon \Omega \to \Rt$ and magnetizations $\boldsymbol{m}\colon \Omega \to \St$ which attain only the easy axes (up to sign). Thus, each limiting magnetization determines a domain structure. 

Theorem \ref{thm:gamma} provides the $\Gamma$-limit $\mathcal{E}$ of the sequence $(\mathcal{E}_\varepsilon)_{\varepsilon>0}$, as $\varepsilon \to 0^+$.  The functional $\mathcal{E}$ is the sum of two contributions (see \eqref{eq:E}--\eqref{eq:Em} below): a linearized elastic energy where the strain corresponding to stress-free configuration is quadratic in the magnetization;  and an energy measuring the area of the interfaces between different magnetic domains  that  resembles the classical $\Gamma$-limit of the vectorial Modica-Mortola functional \cite{baldo}.  The elastic contribution coincides with the one commonly adopted in small-strain models (see, e.g., \cite{james.desimone}), while the magnetic one agrees with the sharp-interface limit identified in \cite{anzellotti.baldo.visintin}. \EEE  

From Theorem \ref{thm:compactness} and Theorem \ref{thm:gamma}, we immediately deduce the convergence of  almost \EEE minimizers by standard $\Gamma$-convergence arguments (see, e.g, \cite[Theorem 7.8]{dalmaso}). As the contributions given by stray and external fields are both continuous with respect to the relevant notion of convergence, these are easily  incorporated \EEE (see, e.g., \cite[Remark 1.7]{braides}). The   convergence   result  for  almost minimizers is given in Corollary \ref{cor:am}.

Eventually, Theorem \ref{thm:regularity} addresses the regularity of minimizers of the functional $\mathcal{E}$ in the case of uniaxial materials  as in \eqref{eq:uniaxial}. \EEE  For sufficiently regular boundary datum and  external field, we prove the regularity of the infinitesimal displacement and the boundaries of magnetic domains. \MMM This result \EEE provides a slight   extension of \cite[Theorem 5.1]{anzellotti.baldo.visintin}  which accounts \EEE for magnetostriction in the small-strain setting.

\subsection{Comments  and outlook} \label{subsec:discussion}
Our   main findings represent the first results combining linearization and sharp-interface limit for Eulerian-Lagrangian energies, i.e., energies comprising both terms in the reference configuration and in the unknown deformed configuration.

The expression of the spontaneous strain of magnetizations postulated here (see \eqref{eq:Lambda} below) is evidenced in the literature (see, e.g.,  \cite[Eq.~(90)]{james.kinderlehrer}) and seems more natural  compared with  the one in \cite{almi.kruzik.molchanova}, even though the two are asymptotically equivalent. In contrast with \cite{gkms-arma,gkms}, the anisotropy energy depends on the deformation gradient,
 so that our diffuse-interface model complies with frame indifference. In this regard, the rigidity of the model arising from the linearization process is steadily exploited while handling this additional dependence on the deformation gradient.

From a technical point of view, we observe that  $p\geq 3$ is assumed in \cite{almi.kruzik.molchanova,gkms-arma,gkms}, while here we only require $p>2$. Thus,  admissible deformation are possibly discontinuous in our setting. The consequent difficulties are overcome by means of the techniques developed in \cite{barchiesi.henao.moracorral,bresciani,bresciani.friedrich.moracorral}. 

Condition \eqref{eq:beta} identifies the regime where no interaction between elastic and magnetic energy term occurs during the process of $\Gamma$-limit. As a result, the magnetic term of the limiting energy \MMM is \EEE independent of the displacement. 
 Different regimes potentially leading to a $\Gamma$-limit where the magnetic term depends also on the displacement will be the subject of future investigations. 
 
 Other interesting and ambitious directions might concern the linearization for spontaneous strains of magnetizations that are not perturbations of the identity, and the sharp-interface limit at large-strains with the anisotropy energy depending on the deformation gradient.

\EEE

\subsection{Structure of the paper} The paper is organized into four sections. Section \ref{sec:intro} constitutes this introduction. In Section \ref{sec:main} we specify the setting of our problem and we state our main findings. Section \ref{sec:preliminaries} collects results from the literature that will be needed for the proofs of our main results presented in Section \ref{sec:proofs}. 

 \EEE 
 
 \subsection{Notation} Given $a,b\in \R$, we employ the notation $a \wedge b\coloneqq \min \{a,b\}$ and $a \vee b\coloneqq \max\{ a,b\}$.
 We work in $\Rt$ with $\St$ denoting the unit sphere centered at the origin. The class $\Rtt$ of square matrices  is endowed with the Frobenius norm. The symbol  $\boldsymbol{id}$ is used for the identity map on $\Rt$, while $\boldsymbol{I}$ and $\boldsymbol{O}$ denote the identity and the null matrix in $\Rtt$. \EEE   
  Given $\boldsymbol{A}\in\Rtt$, we write  $\det \boldsymbol{A}$ for its determinant and  $\boldsymbol{A}^\top$ for its transpose. $GL_+(3)$ is the set of matrices with positive determinant, while  $SO(3)$ is the set of rotations. $\mathrm{Sym}(3)$ stands for the class of symmetric matrices and  $\mathrm{sym}\boldsymbol{A}\coloneqq (\boldsymbol{A}+\boldsymbol{A}^\top)/2$ is the symmetric part of $\boldsymbol{A}\in\Rtt$. We write $E\triangle F$ for the symmetric difference of two sets $E,F\subset \Rt$.
 
 The Lebesgue measure and the two-dimensional Hausdorff measure on $\Rt$ are denoted by $\leb$ and $\haus$, respectively. The expressions ``almost everywhere'' (abbreviated as `a.e.' within the formulas) as well as the adjective `negligible' always refer to the Lebesgue measure. The symbol $\chi_A$ is employed for the indicator function of a set $A\subset \Rt$.
 We use standard notation for Lebsgue spaces ($L^p$), spaces of differentiable and H\"{o}lder functions ($C^k$ and $C^{k,\alpha}$), Sobolev spaces ($W^{k,p}$), and spaces of functions with bounded variation ($BV$). The subscripts `loc' and `c' are used for their local and compactly-supported versions. Domain and codomain are separated by a semicolon and the codomain is omitted for scalar-valued maps.  We write $g\,\d\boldsymbol{x}$ for the measure with density $g$ with respect to $\leb$. \EEE  Lebesgue spaces on boundaries are always meant with respect to $\haus$.  Given $1<p<\infty$, we write $p'\coloneqq p/(p-1)$. \EEE Distributional gradients are indicated   by \EEE $D$, while $E$ is used for symmetrized gradients. \EEE Boundary conditions for Sobolev functions are always understood in the sense of traces. Given a map ${w}$ with bounded variation, we write $|D w|$ and $J_w$ for its total variation and jump set, respectively. Eventually, the class of sets of locally finite perimeter in $\Rt$ is written as $\mathcal{P}(\Rt)$ with $\mathrm{Per}(E;\Omega)$ standing for the relative perimeter of $E\in\mathcal{P}(\Rt)$ in an open set $\Omega\subset \Rt$.
 
  We adopt the usual  convention of denoting by $C$ a positive constant depending on known quantities (in particular, independent on $\varepsilon$) whose value can change from line to line. \EEE 

\section{Setting and main results}

\label{sec:main}

\subsection{Setting}
Throughout the paper, we fix  
\begin{align*}
	 &\hspace{0,5cm}\text{$\Omega\subset \Rt$  bounded domain with Lipschitz boundary,  \quad $p>2$, \quad $q>1$,} \\
		&\hspace{0cm}\text{$\Delta \subset \partial \Omega$  with Lipschitz boundary in $\partial \Omega$}, \quad \haus(\Delta)>0, \quad  \boldsymbol{d}\in W^{1,\infty}(\Omega;\Rt). 	
\end{align*} 
The fact that $\Delta$ has  Lipschitz boundary in $\partial \Omega$ is understood in the sense of \cite[Definition 2.1]{agostiniani.dalmaso.desimone}.
\EEE 

Fix $\beta>0$. For each $\varepsilon>0$, we consider  
the functional
\begin{equation}
	\label{eq:Eeps}
	\mathcal{E}_{\varepsilon}(\boldsymbol{y},\boldsymbol{m})\coloneqq \mathcal{E}_{\varepsilon}^{\rm e}(\boldsymbol{y},\boldsymbol{m})+\mathcal{E}_{\varepsilon}^{\rm m}(\boldsymbol{y},\boldsymbol{m})
\end{equation}
with
\begin{align}
			\label{eq:Ee-eps}
		\mathcal{E}_{\varepsilon}^{\rm e}(\boldsymbol{y},\boldsymbol{m})&\coloneqq \frac{1}{\varepsilon^2}\int_\Omega W(\boldsymbol{\Lambda}_\varepsilon^{-1}(\boldsymbol{m}\circ \boldsymbol{y})D\boldsymbol{y})\,\d\boldsymbol{x},\\
			\label{eq:Em-eps}
		\mathcal{E}_{\varepsilon}^{\rm m}(\boldsymbol{y},\boldsymbol{m})&\coloneqq\frac{1}{\varepsilon^\beta}\int_{\imt(\boldsymbol{y},\Omega)}  \Phi \left( \left( (D\boldsymbol{y})^\top \circ \boldsymbol{y}^{-1}\right ) \boldsymbol{m}   \right)  \,\d\boldsymbol{\xi}+\varepsilon^\beta \int_{\imt(\boldsymbol{y},\Omega)} |D\boldsymbol{m}|^2\,\d\boldsymbol{\xi},
\end{align} \EEE 
defined on the class of admissible states
\begin{equation*}
	\begin{split}
		\mathcal{A}_\varepsilon\coloneqq  \big \{ (\boldsymbol{y},\boldsymbol{m}):\:\boldsymbol{y}\in \mathcal{Y}_{p,q}(\Omega),\:\: \text{$\boldsymbol{y}=\boldsymbol{id}+\varepsilon \boldsymbol{d}$ on $\Delta$,}\:\:\boldsymbol{m}\in W^{1,2}(\imt(\boldsymbol{y},\Omega);\St)&\big \},
	\end{split} 
\end{equation*}
 where
\begin{equation*}
	\mathcal{Y}_{p,q}(\Omega)\coloneqq \left\{ \boldsymbol{y}\in W^{1,p}(\Omega;\Rt):\:\det D\boldsymbol{y}\in L^q(\Omega;(0,+\infty)),\:\:\text{$\boldsymbol{y}$ satisfies (DIV)},\:\:\text{$\boldsymbol{y}$ a.e. injective}   \right\}.
\end{equation*}
The divergence identities (DIV), introduced in \cite{mueller.wcd}, are  given in Definition \ref{def:div} below.  Intuitively, the validity of these identities excludes the phenomenon of cavitation \cite{ball.1982}. \EEE  Recall that  $\boldsymbol{y}$ is termed almost everywhere injective whenever the restriction of $\boldsymbol{y}$ to the complement of a negligible subset of $\Omega$ is an injective map \cite{ball.gi}.  \EEE

The elastic energy  \eqref{eq:Ee-eps} depends on a material density $W\colon \Rtt \to [0,+\infty]$ on which we make the following assumptions:
\begin{enumerate}[label=(W\arabic*)]
	\item \emph{Orientation preservation}: $W(\boldsymbol{A})<+\infty$ if and only if  $ \boldsymbol{A}\in GL_+(3)$;
	\item \emph{Regularity}: $W$ is continuous on $\Rtt$ and of class $C^2$ \EEE in a neighborhood of $SO(3)$;
	\item \emph{Frame indifference}: $W$ is frame indifferent, i.e.,
	\begin{equation*}
		W(\boldsymbol{Q}\boldsymbol{A})=W(\boldsymbol{A}) \quad \text{for all $\boldsymbol{Q}\in SO(3)$ and $ \boldsymbol{A}\in GL_+(3)$;}
	\end{equation*}
	\item  \emph{Ground state}:
	 $W(\boldsymbol{I})=0$;
	\item \emph{Growth}: there exist a constant  $c_W>0$   
	such that
	\begin{equation*}
		\text{$W(\boldsymbol{A})\geq  c_W \EEE  \, \left( g_p(\dist(\boldsymbol{A};SO(3)))+h_q(\det \boldsymbol{A})   \right) $ \quad for all $ \boldsymbol{A}\in GL_+(3)$,}
	\end{equation*}
	where $g_p\colon [0,+\infty) \to[0,+\infty)$ and $h_q\colon (0,+\infty)\to  [  0,+\infty)$ are defined as 
	\begin{equation*}
		g_p(t)\coloneqq \begin{cases}
			\frac{t^2}{2} & \text{if $0\leq t \leq 1$}, \\
			\frac{t^p}{p}+\frac{1}{2}-\frac{1}{ p } & \text{if $t>1$},
		\end{cases}
	\end{equation*}
	and
	\begin{equation*}
		h_q(t)\coloneqq \frac{t^q}{q}+\frac{1}{t}-\frac{q+1}{q} \quad \text{for all $t>0$.}
	\end{equation*}
\end{enumerate}
The specific expressions of the functions $g_p$ and $h_q$ are not relevant as long as $g_p(t)$ is bounded form below by  $t^2 \vee t^p$, and $h_q(t)$ goes as  $t^q$ for $t\gg 1$ and controls ${1}/{t}$ for $t\ll 1$. 
The tensor field $\boldsymbol{\Lambda}_\varepsilon\colon \St \to \Rtt$ appearing in \eqref{eq:Ee-eps} is defined as
\begin{equation}
	\label{eq:Lambda}
	\boldsymbol{\Lambda}_\varepsilon(\boldsymbol{z})\coloneqq a_\varepsilon \boldsymbol{z}\otimes \boldsymbol{z}+b_\varepsilon\left( \boldsymbol{I}-\boldsymbol{z}\otimes \boldsymbol{z} \right) \quad \text{for all $\boldsymbol{z}\in \St$}
\end{equation}
with $	a_\varepsilon\coloneqq 1+\varepsilon a$ and $b_\varepsilon\coloneqq 1 + \varepsilon b$,
where   $a,b\in\R$. \EEE This  tensor   is invertible and its inverse     reads  as  
\begin{equation}
	\label{eq:Lambda-inverse}
	\boldsymbol{\Lambda}_\varepsilon^{-1}(\boldsymbol{z})=\frac{1}{a_\varepsilon} \boldsymbol{z}\otimes \boldsymbol{z}+\frac{1}{b_\varepsilon} \left(\boldsymbol{I}-\boldsymbol{z}\otimes \boldsymbol{z}\right) \quad \text{for all $\boldsymbol{z}\in \St$}.
\end{equation}  
The magnetic energy \eqref{eq:Em-eps} is given by the sum of two integrals on the topological image $\imt(\boldsymbol{y},\Omega)$ of $\Omega$ under $\boldsymbol{y}$. This set was introduced in \cite{barchiesi.henao.moracorral} and it is here recalled in Definition \ref{def:topim} below. \EEE
 The first integral  features the  function 
$\Phi\colon \Rt \to [0,+\infty)$ which is assumed to satisfy:
\begin{enumerate}[label=($\Phi$\arabic*)]
	\item \emph{Multiple wells}: There exist $M\in\N$ with $M\geq 2$ and  distinct elements  $\boldsymbol{b}_1,\dots,\boldsymbol{b}_M\in\St$ \EEE such that $\{ \boldsymbol{z}\in\St:\:\Phi(\boldsymbol{z})=0  \}= \{ \boldsymbol{b}_1,\dots,\boldsymbol{b}_M  \}\EEE $;
	\item \emph{Regularity}: $\Phi$ is continuous on $\Rt$ and   of class $C^1$ \EEE in a neighborhood of $\St$. 
\end{enumerate}
 From the modeling  point of view, in ($\Phi$1) we should take  $M=2N$ with $N\in\N$,  $  \boldsymbol{b}_{i}  \coloneqq \boldsymbol{a}_i$ for  $i=1,\dots,N$, and $\boldsymbol{b}_i\coloneqq -\boldsymbol{a}_{i-N}$ for $i=N+1,\dots,M$ \EEE where  $\boldsymbol{a}_1,\dots,\boldsymbol{a}_N\in\St$ \EEE denote the easy axes. However, this fact is not relevant for the present analysis. 

By passing, we also remark that the energy in \eqref{eq:Eeps}--\eqref{eq:Em-eps} is frame indifferent in view of (W3)  and the objectivity of $\boldsymbol{\Lambda}_\varepsilon^{-1}$, that is, 
\begin{equation*}
	\label{eq:OBJ}
	\boldsymbol{\Lambda}_\varepsilon^{-1}(\boldsymbol{Q}\boldsymbol{z})=\boldsymbol{Q}\boldsymbol{\Lambda}_\varepsilon^{-1}(\boldsymbol{z})\boldsymbol{Q}^\top \quad \text{for all $\boldsymbol{Q}\in SO(3)$ and $\boldsymbol{z}\in\St$.}
\end{equation*}
  \EEE 

The limit functional 
\begin{equation}
	\label{eq:E}
	\mathcal{E}(\boldsymbol{u},\boldsymbol{m})\coloneqq \mathcal{E}^{\rm e}(\boldsymbol{u},\boldsymbol{m})+\mathcal{E}^{\rm m}(\boldsymbol{u},\boldsymbol{m})
\end{equation}
with
\begin{align}
	\label{eq:Ee}
		\mathcal{E}^{\rm e}(\boldsymbol{u},\boldsymbol{m})&\coloneqq \frac{1}{2}\int_\Omega  Q_W\EEE\left( E\boldsymbol{u}-\boldsymbol{\Lambda}(\boldsymbol{m})  \right)\,\d\boldsymbol{x},\\
		\label{eq:Em}
		\mathcal{E}^{\rm m}(\boldsymbol{m})\EEE &\coloneqq  \frac{1}{2}\EEE  \sum_{ \substack{i,j=1,\dots,M\\  i\neq j }
		 }  \sigma_\Phi^{i,j} \EEE \haus \left( \partial^* \{ \boldsymbol{m}= \boldsymbol{b}_i  \} \cap \partial^*\{ \boldsymbol{m}=  \boldsymbol{b}_j  \}   \right),
\end{align}
 is defined on the class
\begin{equation*}
	\begin{split}
		\mathcal{A}\coloneqq \big \{ (\boldsymbol{u},\boldsymbol{m}):\:\:\boldsymbol{u}\in W^{1,2}(\Omega;\Rt),\:\:\text{$\boldsymbol{u}= \boldsymbol{d}\EEE$ on $\Delta$,}\:\:\boldsymbol{m}\in  BV(\Omega;  \{\boldsymbol{b}_1,\dots,\boldsymbol{b}_M  \}  ) 
		\big \}
	\end{split}
\end{equation*}
In \eqref{eq:Ee},   the quadratic form  $Q_W\colon \Rtt \to [0,+\infty)$ is defined as 
\begin{equation*}
	\label{eq:CW}
	Q_W( \boldsymbol{B}\EEE)\coloneqq \mathbf{C}_W( \boldsymbol{B}\EEE): \boldsymbol{B}\EEE \quad \text{for all $\boldsymbol{B}\EEE\in\RNN$,}
\end{equation*} 
where  the elasticity tensor \EEE  $\mathbf{C}_W\colon \RNN \to \RNN$ is defined as the differential of $W$ at the identity, \EEE 
$E\boldsymbol{u}\coloneqq \sym\,D\boldsymbol{u}$ is the symmetrized gradient of $\boldsymbol{u}$,  and, similarly to \eqref{eq:Lambda},  we define $\boldsymbol{\Lambda}\colon \St \to \Rtt$ by setting
\begin{equation}
	\label{eq:Lambdal}
	\boldsymbol{\Lambda}(\boldsymbol{z})\coloneqq a\boldsymbol{z}\otimes \boldsymbol{z}+b\left( \boldsymbol{I}-\boldsymbol{z}\otimes \boldsymbol{z}\right) \quad \text{for all $\boldsymbol{z}\in \St$}.
\end{equation}
Moreover, we set
\begin{equation*}
	\sigma^{i,j}_\Phi\EEE \coloneqq d_\Phi( \boldsymbol{b}_i,\boldsymbol{b}_j\EEE) \quad \text{for all $i,j=1,\dots,M$,}
\end{equation*}
where $d_\Phi\colon \St \times \St \to [0,+\infty)$ denotes the  geodesic \EEE  distance defined   by 
 \begin{equation}
 	\label{eq:distance}
 	d_\Phi( \boldsymbol{z}_0,\boldsymbol{z}_1)\coloneqq \inf \left  \{ \int_0^1 \sqrt{\Phi(\boldsymbol{\gamma}(t))}|\boldsymbol{\gamma}'(t)|\,\d t:\:\:\boldsymbol{\gamma}\in C^1([0,1];\St),\:\:\boldsymbol{\gamma}(0)= \boldsymbol{z}_0,\:\:\boldsymbol{\gamma}(1)= \boldsymbol{z}_1 \right  \}
 \end{equation} 
 for all $ \boldsymbol{z}_0,\boldsymbol{z}_1\in\St$.

\subsection{Main results}
Our main results show that the asymptotic behavior of the functionals $(\mathcal{E}_\varepsilon)_{\varepsilon >0 }$, as $\varepsilon \to 0^+$, is effectively described by the functional $\mathcal{E}$ in the regime \eqref{eq:beta}.
 In this setting, each of the two functionals in \eqref{eq:Ee}--\eqref{eq:Em} is the $\Gamma$-limit of the corresponding term in \eqref{eq:Ee-eps}--\eqref{eq:Em-eps}.     

The first result ensures that the sequence $(\mathcal{E}_\varepsilon)_{\varepsilon >0}$ enjoys suitable compactness properties.

\begin{theorem}[Compactness]
	\label{thm:compactness}
Assume {\rm (W1)--(W5)}, {\rm ($\Phi$1)--($\Phi$2)}, and \eqref{eq:beta}. Let $((\boldsymbol{y}_\varepsilon,\boldsymbol{m}_\varepsilon))_{\varepsilon >0}$ be a sequence with  $(\boldsymbol{y}_\varepsilon,\boldsymbol{m}_\varepsilon)\in \mathcal{A}_\varepsilon$ for all $\varepsilon >0$. Suppose that
\begin{equation}
	\label{eq:E-bdd}
	\sup_{\varepsilon>0} \mathcal{E}_\varepsilon(\boldsymbol{y}_\varepsilon,\boldsymbol{m}_\varepsilon)<+\infty.
\end{equation}
Then, there  exists  $(\boldsymbol{u},\boldsymbol{m})\in\mathcal{A}$  such that, up to subsequences, we have  as $\varepsilon \to 0^+$
\begin{align}
	\label{eq:y}
	\boldsymbol{y}_\varepsilon \to  \boldsymbol{id} \quad &\text{in $W^{1,p}(\Omega;\Rt)$,}\\
	\label{eq:u}
	 \frac{1}{\varepsilon} \left( \boldsymbol{y}_\varepsilon - \boldsymbol{id} \right) \wk \boldsymbol{u} \quad &\text{in $W^{1,2}(\Omega;\Rt)$},\\
	\label{eq:m}
	\chi_{\imt(\boldsymbol{y}_\varepsilon,\Omega)}\boldsymbol{m}_\varepsilon \to \chi_\Omega\boldsymbol{m} \quad &\text{in $L^1(\Rt;\Rt)$,}\\
		\label{eq:composition}
	\text{$\boldsymbol{m}_\varepsilon \circ \boldsymbol{y}_\varepsilon \to \boldsymbol{m}$} \quad &\text{in $L^1(\Omega;\Rt)$.}
\end{align}
\end{theorem}

The second result shows that the functional $\mathcal{E}$ provides an optimal lower bound for $(\mathcal{E}_\varepsilon)_{\varepsilon>0}$. For the notion of $\Gamma$-convergence, we refer to \cite{braides,dalmaso}.
 
\begin{theorem}[$\Gamma$-convergence]
	\label{thm:gamma}
Assume {\rm (W1)--(W5)}, {\rm ($\Phi$1)--($\Phi$2)}, and \eqref{eq:beta}.   Then, the sequence  $(\mathcal{E}_\varepsilon)_{\varepsilon>0}$ $\Gamma$-converges to the functional $\mathcal{E}$, as $\varepsilon\to 0^+$, with respect to the convergences in Theorem \ref{thm:compactness}. Namely, the following two conditions hold:
\begin{enumerate}[label=(\roman*)]
	\item \emph{Lower bound:} For any sequence $\left( (\boldsymbol{y}_\varepsilon,\boldsymbol{m}_\varepsilon) \right)_{\varepsilon>0}$ with $(\boldsymbol{y}_\varepsilon,\boldsymbol{m}_\varepsilon)\in\mathcal{A}_\varepsilon$ for all $\varepsilon>0$ and any $(\boldsymbol{u},\boldsymbol{m})\in\mathcal{A}$ for which \eqref{eq:y}--\eqref{eq:composition} hold  true, we have 
	\begin{equation}
		\label{eq:lb}
		\liminf_{\varepsilon \to 0^+} \mathcal{E}_\varepsilon(\boldsymbol{y}_\varepsilon,\boldsymbol{m}_\varepsilon)\geq \mathcal{E}(\boldsymbol{u},\boldsymbol{m});
	\end{equation} 
	\item \emph{Optimality of the lower bound:} For any $(\boldsymbol{u},\boldsymbol{m})\in\mathcal{A}$, there exists a sequence $\left( (\boldsymbol{y}_\varepsilon,\boldsymbol{m}_\varepsilon) \right)_{\varepsilon>0}$ with $(\boldsymbol{y}_\varepsilon,\boldsymbol{m}_\varepsilon)\in\mathcal{A}_\varepsilon$ for all $\varepsilon>0$ for which  \eqref{eq:y}--\eqref{eq:composition} hold true and we have
	\begin{equation}
		\label{eq:olb}
		\lim_{\varepsilon \to 0^+} \mathcal{E}_\varepsilon(\boldsymbol{y}_\varepsilon,\boldsymbol{m}_\varepsilon)= \mathcal{E}(\boldsymbol{u},\boldsymbol{m}).
	\end{equation}
\end{enumerate}
\end{theorem}

 The effect of magnetic fields can be easily incorporated within our analysis.    Following \cite{desimone.podioguidugli}, we define the stray field corresponding to a deformation $\boldsymbol{y}\in \mathcal{Y}_{p,q}(\Omega)$  and a magnetization  $\boldsymbol{m}\in L^1(\imt(\boldsymbol{y},\Omega);\St)$ as the unique  distributional  solution $\boldsymbol{h}^{\boldsymbol{y},\boldsymbol{m}}\in L^2(\Rt;\Rt)$  to \EEE the  Maxwell system
\begin{equation}
	\label{eq:maxwell}
	\begin{cases}
		\mathrm{curl}\,\boldsymbol{h}^{\boldsymbol{y},\boldsymbol{m}}=\boldsymbol{0} & {}\\
		\mathrm{div}\,\boldsymbol{h}^{\boldsymbol{y},\boldsymbol{m}}=-\mathrm{div}\left( \chi_{\imt(\boldsymbol{y},\Omega)}  \boldsymbol{m} \det D \boldsymbol{y}^{-1} \right)
	\end{cases} \quad \text{in $\Rt$},
\end{equation}
 where $\boldsymbol{y}^{-1}\colon \imt(\boldsymbol{y},\Omega ) \to \Rt$ denotes the inverse of $\boldsymbol{y}$ as given in Proposition \ref{prop:regularity-inverse} below. As a consequence of  Proposition \ref{prop:maxwell} below,  
existence and uniqueness of the stray field are guaranteed in our setting and the associated  energy contribution   
\begin{equation*}
	\mathcal{H}(\boldsymbol{y},\boldsymbol{m})\coloneqq \int_{\Rt} |\boldsymbol{h}^{\boldsymbol{y},\boldsymbol{m}}|^2\,\d\boldsymbol{\xi}
\end{equation*}
is well defined.
Additionally, the work of the  applied magnetic field $\boldsymbol{f}\in L^1(\Rt;\Rt)$ is accounted by the Zeeman energy functional
\begin{equation*}
	\mathcal{F}(\boldsymbol{y},\boldsymbol{m})\coloneqq \int_{\imt(\boldsymbol{y},\Omega)} \boldsymbol{f}\cdot \boldsymbol{m}\,\d\boldsymbol{\xi}.
\end{equation*}
 For  a fixed constant $\lambda\geq 0$, we define the total energies 
\begin{equation*}
	\mathcal{G}_\varepsilon(\boldsymbol{y},\boldsymbol{m})\coloneqq \mathcal{E}_\varepsilon(\boldsymbol{y},\boldsymbol{m}) + \lambda \mathcal{H}(\boldsymbol{y},\boldsymbol{m})-\mathcal{F}(\boldsymbol{y},\boldsymbol{m}) \quad \text{for all $\varepsilon>0$ and $(\boldsymbol{y},\boldsymbol{m})\in\mathcal{A}_\varepsilon$}
\end{equation*}
and
\begin{equation*}
	\mathcal{G}(\boldsymbol{u},\boldsymbol{m})\coloneqq \mathcal{E}(\boldsymbol{u},\boldsymbol{m}) + \lambda \mathcal{H}(\boldsymbol{id},\boldsymbol{m})-\mathcal{F}(\boldsymbol{id},\boldsymbol{m}) \quad \text{for all $(\boldsymbol{u},\boldsymbol{m})\in\mathcal{A}$.}
\end{equation*}

It turns out that \EEE the functionals $\mathcal{H}$ and $\mathcal{F}$  constitute continuous perturbations with respect to relevant topology. \EEE  
Therefore, standard $\Gamma$-convergence arguments   (see, e.g., \cite[Remark 1.7]{braides})  yield the following statement for sequences of almost minimizers.

\begin{corollary}[Convergence of almost minimizers]
	\label{cor:am}
Assume {\rm (W1)--(W5)}, {\rm ($\Phi$1)--($\Phi$2)}, and \eqref{eq:beta}. 
 Let $\left ((\boldsymbol{y}_\varepsilon,\boldsymbol{m}_\varepsilon)\right)_{\varepsilon>0}$ be a sequence with $(\boldsymbol{y}_\varepsilon,\boldsymbol{m}_\varepsilon)\in\mathcal{A}_\varepsilon$ for all $\varepsilon>0$ satisfying
\begin{equation}
	\label{eq:am}
	\lim_{\varepsilon \to 0^+} \left( \mathcal{G}_\varepsilon(\boldsymbol{y}_\varepsilon,\boldsymbol{m}_\varepsilon) - G_\varepsilon  \right)=0,
\end{equation}	
where we set $G_\varepsilon \coloneqq \inf\{ \mathcal{G}_\varepsilon( {\boldsymbol{w}},{\boldsymbol{n}}\EEE):\:( {\boldsymbol{w}},{\boldsymbol{n}}\EEE)\in\mathcal{A}_\varepsilon \} $. 
Then, there exists $(\boldsymbol{u},\boldsymbol{m})\in\mathcal{A}$  such that, up to subsequences, \eqref{eq:y}--\eqref{eq:composition} hold true  and  $(\boldsymbol{u},\boldsymbol{m})$ is a minimizer of  $\mathcal{G}$   with  
\begin{equation}
	\label{eq:GGG}
	\lim_{\varepsilon \to 0^+} \mathcal{G}_\varepsilon(\boldsymbol{y}_\varepsilon,\boldsymbol{m}_\varepsilon)=\mathcal{G}(\boldsymbol{u},\boldsymbol{m}).
\end{equation}
\end{corollary}

 Eventually, we discuss the regularity of minimizers for the limiting energy in a special case. For simplicity, as in \cite{anzellotti.baldo.visintin},  we discuss only the interior regularity in the case $M=2$. This setting can describe uniaxial materials with $\Phi$ as in \eqref{eq:uniaxial} by choosing $\boldsymbol{b}_1\coloneqq  \boldsymbol{a}$ and $\boldsymbol{b}_2\coloneqq -\boldsymbol{a}$.
 In this situation, the limiting magnetic energy \eqref{eq:Em}  can be trivially rewritten as
 \begin{equation}
 	\label{eq:due}
 	\mathcal{E}^{\rm m}(\boldsymbol{m})=\MMM {\sigma_\Phi^{1,2}} \EEE \,  \per (\{ \boldsymbol{m}=\boldsymbol{b}_1 \};\Omega).
 \end{equation} 
Therefore,  we can resort to the regularity theory for sets with minimal perimeter \cite{tamanini}.
  \EEE 

\begin{theorem}[Regularity of minimizers]
	\label{thm:regularity}
Let $( \boldsymbol{u},\boldsymbol{m}\EEE)\in\mathcal{A}$ be a minimizer of the functional $\mathcal{G}$. Then, 	$ \boldsymbol{u}\EEE\in W^{1,r}(\Omega;\Rt)$ for all $r\in [1,\infty)$. Additionally, assume $M=2$. If $\boldsymbol{f}\in L^{s_{\mathrm{f}}} (\Rt;\Rt)$ with $s_{\rm f}> 3$, then  the set $J_{ \boldsymbol{m} \EEE}$ is  a surface \EEE of class $C^{1,(s_{\rm f}-3)/(2s_{\rm f})}$ in $\Omega$. Furthermore,  if   $\boldsymbol{d}\in W^{2,s_{\rm d}}_{\rm loc}(\Omega;\Rt)$ with $s_{\rm d}>3$, then $ \boldsymbol{u} \EEE\in C^{1,1-{3}/{s_{\rm d}}}(\Omega\setminus J_{ \boldsymbol{m}\EEE};\Rt)$. 
\end{theorem}

\section{Preliminary results}

\label{sec:preliminaries}

In this section, we  collect some preliminary notions and facts that will be instrumental for the proofs of the results in Section \ref{sec:main}.

\subsection{Geometric and topological image} 
We introduce two concepts of image for deformations, namely that of geometric  and topological image. For the notions of   approximate differentiability and density of a set, we refer to \cite{ambrosio.fusco.pallara}. \EEE

\begin{definition}[Geometric image]
Let $\boldsymbol{y}\in W^{1,p}(\Omega;\Rt)$ with $\det D \boldsymbol{y}(\boldsymbol{x})>0$ for almost all $\boldsymbol{x}\in\Omega$. The geometric domain $\domg(\boldsymbol{y},\Omega)$ of $\boldsymbol{y}$ is defined as the set of points $\boldsymbol{x}\in \Omega$ at which $\boldsymbol{y}$ is approximately differentiable  with $\det D \boldsymbol{y}(\boldsymbol{x})>0$	and there exist a compact set $K\subset \Omega$ with density one at $\boldsymbol{x}$ and a map $\boldsymbol{w}\in C^1(\Rt;\Rt)$ such that $\boldsymbol{w}\restr{K}=\boldsymbol{y}\restr{K}$ and $(D\boldsymbol{w})\restr{K}=(D\boldsymbol{y})\restr{K}$. Moreover, for any measurable set $A\subset \Omega$, the geometric image of $A$ under $\boldsymbol{y}$ is defined as $\img(\boldsymbol{y},A)\coloneqq \boldsymbol{y}\left( A\cap \domg(\boldsymbol{y},\Omega)  \right)$.
\end{definition}

From \cite[p.~582--583]{henao.moracorral.fracture}, we known that $\domg(\boldsymbol{y},\Omega)$ has full measure in $\Omega$, while \cite[Lemma 3]{henao.moracorral.fracture} ensures that $\boldsymbol{y}\restr{\domg(\boldsymbol{y},\Omega)}$ is (everywhere) injective. Hence, we can consider the inverse of $\boldsymbol{y}$ as the map $\boldsymbol{y}^{-1}\colon \img(\boldsymbol{y},\Omega)\to \domg(\boldsymbol{y},\Omega)$.

Below, $\boldsymbol{y}^*\colon \Omega \to \Rt$ will denote the precise representative of the map $\boldsymbol{y}\in W^{1,p}(\Omega;\Rt)$ defined as
\begin{equation*}
	\boldsymbol{y}^*(\boldsymbol{x})\coloneqq \limsup_{r\to 0^+} \dashint_{B(\boldsymbol{x},r)} \boldsymbol{y}({\boldsymbol{\chi}})\,\d{\boldsymbol{\chi}}\quad \text{for all $\boldsymbol{x}\in\Omega$.}
\end{equation*}
In loose terms,  the restriction of $\boldsymbol{y}^*$ to the boundary of almost every set $U\subset \subset \Omega$ is continuous as a consequence of the assumption $p>2$ and Morrey's embedding. This observation allows for a notion of topological degree for Sobolev mappings. For the  degree of continuous functions, we refer to \cite{fonseca.gangbo}. The next  definition  exploits the sole dependence of the degree on boundary values.

\begin{definition}[Topological image]
	\label{def:topim}
	Let $\boldsymbol{y}\in W^{1,p}(\Omega;\Rt)$ and  $U\subset \subset \Omega$ be a domain such that $\boldsymbol{y}^*\restr{\partial U}\in C^0(\partial U;\Rt)$. The topological image of $U$ under $\boldsymbol{y}$   is defined   as
	\begin{equation*}
		\imt(\boldsymbol{y},U)\coloneqq \left\{ \boldsymbol{\xi}\in\Rt \setminus {\boldsymbol{y}}^*(\partial U):\:\deg({\boldsymbol{y}}^*,U,\boldsymbol{\xi})\neq 0  \right\},
	\end{equation*}
	where  $\deg({\boldsymbol{y}}^*,U,\boldsymbol{\xi})$ denotes  the topological degree of any extension of ${\boldsymbol{y}}^*\restr{\partial U}$ to $\closure{U}$ at the point $\boldsymbol{\xi}$.  
	Moreover, the topological image of $\Omega$  under $\boldsymbol{y}$ is defined as
	\begin{equation*}
		\imt(\boldsymbol{y},\Omega)\coloneqq \bigcup \left\{ \imt(\boldsymbol{y},U):\:\text{$U\subset \subset \Omega$  domain with $\boldsymbol{y}^*\restr{\partial U}\in C^0(\partial U;\Rt)$} \right\}.
	\end{equation*}
\end{definition}
By classical properties of the topological degree (see, e.g.,  \cite[Remark 2.17]{bresciani.friedrich.moracorral}), $\imt(\boldsymbol{y},U)$ is a bounded open set. The set $\imt(\boldsymbol{y},\Omega)$ is thus open but possibly unbounded. 

\subsection{Divergence identities}

The divergence identities are a generalization of the  identity $\Det D \boldsymbol{y}=\det D \boldsymbol{y}$, asserting   that the distributional and pointwise determinant coincide \cite{mueller.det}, which    have been investigated in connection with the weak continuity of the Jacobian determinant \cite{mueller.wcd}. 

\begin{definition}[Divergence identities]
	\label{def:div}
Let $\boldsymbol{y}\in W^{1,p}(\Omega;\Rt)$. The map $\boldsymbol{y}$ is termed to satisfy the divergence identities {\rm (DIV)} whenever
\begin{equation*}
	\div \left( (\adj D\boldsymbol{y})\boldsymbol{\psi}\circ \boldsymbol{y}  \right)=(\div \boldsymbol{\psi})\circ \boldsymbol{y}\det D \boldsymbol{y}\quad \text{for all $\boldsymbol{\psi}\in C^1_{\rm c}(\Rt;\Rt)$,}
\end{equation*}
where the divergence on the left-hand side and the equality are understood in the sense of distributions over $\Omega$. Namely,
\begin{equation*}
	-\int_\Omega \left( (\adj D \boldsymbol{y})\boldsymbol{\psi}\circ \boldsymbol{y} \right)\cdot D\varphi \,\d\boldsymbol{x}=\int_\Omega (\div \boldsymbol{\psi})\circ \boldsymbol{y} (\det D \boldsymbol{y})\varphi\,\d\boldsymbol{x} 
\end{equation*}
for all $\varphi \in C^1_{\rm c}(\Omega)$ and $\boldsymbol{\psi}\in C^1_{\rm c}(\Rt;\Rt)$.
\end{definition}
Let us  remark  that satisfying the divergence identities is equivalent to having zero surface energy in the setting of \cite{barchiesi.henao.moracorral,henao.moracorral.invertibility}.

We state two consequences of the divergence identities.
The first result  shows that geometric and topological image actually coincide and, in turn, rules out the phenomenon of cavitation. \EEE 

\begin{proposition}[Excluding cavitation]
	\label{prop:img=imt}
Let $\boldsymbol{y}\in W^{1,p}(\Omega;\Rt)$ with $\det D \boldsymbol{y}\in L^1(\Omega;(0,+\infty))$ satisfy {\rm (DIV)}. Then,  $\img(\boldsymbol{y},\Omega)=\imt(\boldsymbol{y},\Omega)$ up to negligible sets. That is, $\leb \left( \img(\boldsymbol{y},\Omega)\triangle \imt(\boldsymbol{y},\Omega) \right)=0$.
\end{proposition}
\begin{proof}
See \cite[Lemma 5.18(c)]{barchiesi.henao.moracorral}.
\end{proof}

The second result concerns the inverse of deformations. 
As already noted, if $\boldsymbol{y}$ is almost everywhere injective, then the restriction $\boldsymbol{y}\restr{\domg(\boldsymbol{y},\Omega)}$ is actually (everywhere) injective. Thus, the inverse  $\boldsymbol{y}^{-1}\colon \imt(\boldsymbol{y},\Omega)\to \domg(\boldsymbol{y},\Omega)$   yields an almost everywhere defined map on the open set $\imt(\boldsymbol{y},\Omega)$ in view of Proposition \ref{prop:img=imt}. The next results asserts its Sobolev  regularity. 

\begin{proposition}[Regularity of the inverse]
	\label{prop:regularity-inverse}
Let $\boldsymbol{y}\in W^{1,p}(\Omega;\Rt)$ with $\det D \boldsymbol{y}\in L^1(\Omega;(0,+\infty))$ be almost everywhere injective and satisfy {\rm (DIV)}. Then, $\boldsymbol{y}^{-1}\in W^{1,1}(\imt(\boldsymbol{y},\Omega);\Rt)$ with $D\boldsymbol{y}^{-1}=(D\boldsymbol{y})^{-1}\circ \boldsymbol{y}^{-1}$ almost everywhere in $\imt(\boldsymbol{y},\Omega)$.
\end{proposition}
\begin{proof}
	From \cite[Proposition 5.3]{barchiesi.henao.moracorral}, we know that $\boldsymbol{y}^{-1}\in W^{1,1}_{\rm loc}(\imt(\boldsymbol{y},\Omega);\Rt)$ and the inverse formula holds. Moreover,
	\begin{equation*}
		\int_{\imt(\boldsymbol{y},\Omega)} |\boldsymbol{y}^{-1}|\,\d\boldsymbol{\xi}=\int_\Omega |\boldsymbol{id}|\det D \boldsymbol{y}\,\d\boldsymbol{x}\leq 	\|\boldsymbol{id}\|_{L^\infty(\Omega;\Rt)}\,\|\det D\boldsymbol{y}\|_{L^1(\Omega)}
	\end{equation*}
	and
	\begin{equation*}
		\int_{\imt(\boldsymbol{y},\Omega)} |D\boldsymbol{y}^{-1}|\,\d\boldsymbol{\xi}=\int_\Omega |\adj D \boldsymbol{y}|\,\d\boldsymbol{x}
	\end{equation*}
	thanks to Proposition \ref{prop:img=imt}, the inverse formula, and the change-of-variable formula,  where the last integral is finite since $p>2$. \EEE 
\end{proof}

\EEE

\subsection{Least upper bound of positive measures}

 Let $U\subset \Rt$ be an open set and let $(\mu_1,\dots,\mu_M)$ be a finite family  of  positive Borel measures on $U$. The least upper bound $\bigvee_{i=1}^M \mu_i$ of $(\mu_1,\dots,\mu_M)$ is  defined for every Borel set $E\subset U$ by setting
\begin{equation*}
	\label{eq:lubm}
	\left( \bigvee_{i=1}^M \mu_i \right)(E)\coloneqq \sup \left\{ \sum_{i=1}^M \mu_i(E_i):\:\text{$(E_1,\dots,E_M)$ Borel partition of $E$}  \right\}.
\end{equation*}
This is also a positive Borel measure,  namely,  the smallest positive Borel measure on $U$ with the property that
 $\max \{ \mu_i(E):\:i=1,\dots,M \}\leq (\bigvee_{i=1}^M\mu_i)(E)$ for any Borel set $E\subset U$  (see, e.g., \cite[p.~31]{ambrosio.fusco.pallara}). In particular, as noted in \cite[Remark 1.69]{ambrosio.fusco.pallara}, for absolutately continuous measures one has 
\begin{equation}
	\label{eq:vee}
	\bigvee_{i=1}^M (g_i\,\d\boldsymbol{x})=\max\{ g_i:\:i=1,\dots,M\}\,\d\boldsymbol{x} \quad \text{for all $g_1,\dots,g_M\in L^1(U;(0,+\infty))$.}
\end{equation}

\EEE

The following lower semicontinuity result will be employed in the proof of Theorem \ref{thm:gamma}(i).  

\EEE 

\begin{lemma}[Lower semicontinuity of the least upper bound of total variations]
	\label{lem:tot-var}
	Let $U\subset \Rt$ be an open set and,
	for each $i=1,\dots,M$, let $(w_{i,\varepsilon})_{\varepsilon >0}$ be a sequence in $BV(U)$ and $w_i\in BV(U)$ such that $w_{i,\varepsilon}\to w_i$ in $L^1_{\rm loc}(U)$, as $\varepsilon \to 0^+$. Then
	\begin{equation*}
		\left( \bigvee_{i=1}^M |Dw_i| \right)(U)\leq \liminf_{\varepsilon \to 0^+} \left( \bigvee_{i=1}^M |Dw_{i,\varepsilon}| \right)(U).
	\end{equation*} 	
\end{lemma}
\begin{proof}
	Let  $(E_1,\dots,E_M)$ be a Borel  partition of $U$. By the lower semicontinuity of the total variation,	
	\begin{equation*}
		|Dw_i|( E_i\EEE)\leq \liminf_{\varepsilon \to 0^+} |Dw_{i,\varepsilon}|( E_i\EEE)\quad \text{for all $i=1,\dots,M$.}
	\end{equation*}
Thus
	\begin{equation*}
		\begin{split}
			\sum_{i=1}^M |Dw_i|( E_i\EEE)&\leq \sum_{i=1}^M \liminf_{\varepsilon \to 0^+} |Dw_{i,\varepsilon}|( E_i\EEE)\\
			&\leq \liminf_{\varepsilon \to 0^+} \left(  \sum_{i=1}^M |Dw_{i,\varepsilon}|( E_i\EEE) \right)\leq \lim_{\varepsilon \to 0^+} \left( \bigvee_{i=1}^M |Dw_{i,\varepsilon}| \right)(U).
		\end{split}
	\end{equation*}
	The claim follows by taking the supremum among all Borel partitions $(E_1,\dots,E_M)$ of $U$.
\end{proof}

\section{Proofs}
\label{sec:proofs}

In this section, we prove the main results stated in Section \ref{sec:main}.

\subsection{Compactness result}

We begin by proving our compactness result.  Below, we will consider the first-order Taylor expansion of $\Phi$. For fixed $\boldsymbol{z}_0\in\St$, we have
\begin{equation}
	\label{eq:taylor-Phi}
	\Phi(\boldsymbol{z}_0+\boldsymbol{z})=\Phi(\boldsymbol{w})+D\Phi(\boldsymbol{z}_0)\cdot \boldsymbol{z}+\eta_\Phi(\boldsymbol{z};\boldsymbol{z}_0) \quad \text{for all $\boldsymbol{z}\in \Rt$ with $|\boldsymbol{z}|\ll 1$,}
\end{equation}
where $\eta_\Phi(\boldsymbol{z};\boldsymbol{z}_0)=\mathrm{o}(|\boldsymbol{z}|)$, as $|\boldsymbol{z}|\to 0^+$. For $t>0$, we define
\begin{equation*}
	\zeta_\Phi(t)\coloneqq \sup \left\{ \frac{|\eta_\Phi(\boldsymbol{z};\boldsymbol{z}_0)|}{|\boldsymbol{z}|}:\:\:\boldsymbol{z}_0\in\St,\:\:\boldsymbol{z}\in\Rt,\:\:|\boldsymbol{z}|\leq t  \right\},
\end{equation*}
so that
\begin{equation}
	\label{eq:zeta-Phi}
	\lim_{t\to 0^+} \zeta_\Phi(t)=0
\end{equation}
in view of ($\Phi$2). Additionally,  for $i=1,\dots,M$ we define  $f_i\colon \St \to [0,+\infty)$ by setting 
\begin{equation}
	\label{eq:fi}
	f_i(\boldsymbol{z})\coloneqq d_\Phi(\boldsymbol{z}, \boldsymbol{b}_i) \quad \text{for all $\boldsymbol{z}\in\St$},
\end{equation}
where $d_\Phi$ is given by  \eqref{eq:distance}.

\EEE 

\begin{proof}[Proof of Theorem \ref{thm:compactness}]
We subdivide the proof into three steps.

\emph{Step 1 (Compactness of deformations and displacements).}	
For convenience, we set $\boldsymbol{F}_\varepsilon\coloneqq D\boldsymbol{y}_\varepsilon$, $\boldsymbol{L}_\varepsilon\coloneqq \boldsymbol{\Lambda}_\varepsilon(\boldsymbol{m}_\varepsilon \circ \boldsymbol{y}_\varepsilon)$, and $\boldsymbol{A}_\varepsilon\coloneqq \boldsymbol{L}_\varepsilon^{-1}\boldsymbol{F}_\varepsilon$,
where we adopt the notation in \eqref{eq:Lambda}.
From \eqref{eq:E-bdd} and assumption  (W5), we obtain
\begin{equation}
	\label{eq:E-bdd1}
	 \int_\Omega \left (  g_p(\dist(\boldsymbol{A}_\varepsilon;SO(3)))+h_q(\det \boldsymbol{A}_\varepsilon)  \right ) \,\d \boldsymbol{x}   \leq \frac{\varepsilon^2}{c_W}  \mathcal{E}^{\rm e}_\varepsilon(\boldsymbol{y}_\varepsilon,\boldsymbol{m}_\varepsilon)\leq C \varepsilon^2.  
\end{equation}	
Given that  $t^p/p\leq g_p(t)$  for all $t\geq 0$,
the previous estimate implies $\dist(\boldsymbol{A}_\varepsilon;SO(3))\to 0$ in $ L^p(\Omega)$. Thus,  the sequence $(\boldsymbol{A}_\varepsilon)_{\varepsilon>0}$ is bounded in $ L^p(\Omega;\Rtt)$.   Noting that $\boldsymbol{\Lambda}_\varepsilon(\boldsymbol{z})=\boldsymbol{I}+\varepsilon\boldsymbol{\Lambda}(\boldsymbol{z})$ for all $\boldsymbol{z}\in \St$, where we employ \eqref{eq:Lambda} and \eqref{eq:Lambdal}, and computing \EEE
\begin{equation}
	\label{eq:ll} \text{$|\boldsymbol{\Lambda}(\boldsymbol{z})|=\ell\coloneqq \sqrt{ a^2+2b^2}$  \quad for all $\boldsymbol{z}\in \St$,}	
\end{equation}
the triangle inequality gives
\begin{equation}
	\label{eq:E-bdd2}
	\begin{split}
		\dist(\boldsymbol{F}_\varepsilon;SO(3))&=\dist(\boldsymbol{L}_\varepsilon\boldsymbol{A}_\varepsilon;SO(3))\\
		&\leq |\boldsymbol{L}_\varepsilon -\boldsymbol{I}|\,|\boldsymbol{A}_\varepsilon|+\dist(\boldsymbol{A}_\varepsilon;SO(3))\\
		&= \ell\varepsilon |\boldsymbol{A}_\varepsilon|+\dist(\boldsymbol{A}_\varepsilon;SO(3)).
	\end{split}
\end{equation}
  The estimates \eqref{eq:E-bdd1} and \eqref{eq:E-bdd2} together with the monotonicity of $g_p$ and the elementary inequality
\begin{equation*}
	g_p(t+\tilde{t})\leq C \left(  t^p + g_p(\tilde{t}) \right) \quad \text{for all $t,\tilde{t}\geq 0$},
\end{equation*}
 yield
\begin{equation*}
	\begin{split}
		\int_\Omega g_p(\dist(\boldsymbol{F}_\varepsilon;SO(3)))\,\d\boldsymbol{x}&\leq C \left(  \ell^p \varepsilon^p  \int_\Omega  |\boldsymbol{A}_\varepsilon|^p\EEE \,d\boldsymbol{x} +\int_\Omega g_p(\dist(\boldsymbol{A}_\varepsilon;SO(3)))\,\d\boldsymbol{x}  \right)\leq C \varepsilon^2
	\end{split}
\end{equation*}
given that $p>2$. As
\begin{equation*}
	\frac{1}{2p}(t^2+t^p)\leq \frac{1}{p} (t^2\vee t^p)\leq g_p(t) \quad \text{for all $t\geq 0$,}
\end{equation*} 
from the last estimate   we obtain
\begin{equation}
	\label{eq:rigidity}
	\int_\Omega \dist^2(\boldsymbol{F}_\varepsilon;SO(3))\,\d\boldsymbol{x}+\int_\Omega \dist^p(\boldsymbol{F}_\varepsilon;SO(3))\leq C \varepsilon^2.
\end{equation}
Let $\boldsymbol{u}_\varepsilon\coloneqq \varepsilon^{-1}(\boldsymbol{y}_\varepsilon-\boldsymbol{id})$ and $\boldsymbol{G}_\varepsilon\coloneqq D\boldsymbol{u}_\varepsilon$. First, we show that $(\boldsymbol{G}_\varepsilon)_{\varepsilon >0}$ is bounded in $L^2(\Omega;\Rtt)$. Then, by applying the  Poincar\'{e} inequality with trace term, we find $\boldsymbol{u}\in W^{1,2}(\Omega;\Rt)$ for which  \eqref{eq:u} holds  along a not relabeled subsequence. In particular, \EEE $\boldsymbol{u}=\boldsymbol{d}$ on $\Delta$ \EEE owing to the weak continuity of the trace operator.
 From \eqref{eq:rigidity}, using the classical rigidity estimate \cite[Theorem 3.1]{fjm}, we identify a sequence $(\boldsymbol{R}_\varepsilon)_{\varepsilon >0}$ in $SO(3)$ satisfying
\begin{equation}
	\label{eq:FR}
	\int_\Omega |\boldsymbol{F}_\varepsilon - \boldsymbol{R}_\varepsilon|^2\,\d\boldsymbol{x}\leq C \varepsilon^2.
\end{equation}
Then, by arguing as in \cite[Proposition 3.4]{dalmaso.negri.percivale}, taking advantage of the boundary condition \EEE $\boldsymbol{u}_\varepsilon=\boldsymbol{d}$ on $\Delta$, \EEE we establish
\begin{equation}
	\label{eq:RI}
	|\boldsymbol{R}_\varepsilon -\boldsymbol{I}|^2\leq C \varepsilon^2 \left( 1+ \int_\Gamma |\boldsymbol{d}|^2\,\d\haus  \right)\leq C \varepsilon^2 
\end{equation}
and, in turn,
\begin{equation*}
	\int_\Omega |\boldsymbol{G}_\varepsilon|^2\,\d\boldsymbol{x}= \frac{1}{\varepsilon^2}\int_\Omega |\boldsymbol{F}_\varepsilon-\boldsymbol{I}|^2\,\d\boldsymbol{x}\leq C.
\end{equation*}
This proves the claim. 

Next, we check that $\varepsilon \boldsymbol{G}_\varepsilon\to \boldsymbol{O}$  in $L^p(\Omega;\Rtt)$. Then,    \eqref{eq:y} follows once again thanks to the Poincar\'{e} inequality. We argue as in \cite[Lemma 3.3]{almi.kruzik.molchanova} and   we estimate
\begin{equation*}
	\begin{split}
		\int_\Omega |\varepsilon \boldsymbol{G}_\varepsilon|^p\,\d\boldsymbol{x}&\leq 2^{p-1} \left( \int_\Omega |\boldsymbol{F}_\varepsilon-\boldsymbol{R}_\varepsilon|^p\,\d\boldsymbol{x} + \int_\Omega |\boldsymbol{R}_\varepsilon -\boldsymbol{I}|^p\,\d\boldsymbol{x}  \right).
	\end{split}
\end{equation*}
For the first integral on the right-hand side, as $p>2$, we have
\begin{equation*}
	\int_\Omega |\boldsymbol{F}_\varepsilon-\boldsymbol{R}_\varepsilon|^p\,\d\boldsymbol{x}=\| \boldsymbol{F}_\varepsilon-\boldsymbol{R}_\varepsilon\|_{L^p(\Omega;\Rtt)}^p\leq C \| \boldsymbol{F}_\varepsilon-\boldsymbol{R}_\varepsilon\|_{L^2(\Omega;\Rtt)}^p \leq C \varepsilon^p
\end{equation*}
 thanks to \eqref{eq:FR}, while for the second one, \eqref{eq:RI} gives
\begin{equation*}
	|\boldsymbol{R}_\varepsilon - \boldsymbol{I}|^p=\left( |\boldsymbol{R}_\varepsilon - \boldsymbol{I}|^2 \right)^{p/2}\leq (C\varepsilon^2)^{p/2}=C \varepsilon^p.
\end{equation*}
Altogether, we find 
\begin{equation*}
	\text{$\|\varepsilon \boldsymbol{G}_\varepsilon\|_{L^p(\Omega;\Rtt)}\leq C \varepsilon$}	
\end{equation*}
 which yields the desired convergence.

Further, we claim that 
\begin{equation}
	\label{eq:det1}
	\det D \boldsymbol{y}_\varepsilon \to 1 \quad  \text{in $L^1(\Omega)$},
\end{equation}
so that, by applying \cite[Proposition 2.23(ii)]{bresciani}, we obtain
\begin{equation}
	\label{eq:imt}
	\chi_{\imt(\boldsymbol{y}_\varepsilon,\Omega)}\to \chi_\Omega \quad \text{in $L^1(\Rt)$.}
\end{equation}
Note that
\begin{equation}
	\label{eq:det-Lambda}
	\det \boldsymbol{\Lambda}_\varepsilon(\boldsymbol{z})=a_\varepsilon b_\varepsilon^2 \eqqcolon c_\varepsilon \quad \text{for all $\boldsymbol{z}\in \St$,}
\end{equation}
where $c_\varepsilon \to 1$. Hence, as $\det \boldsymbol{F}_\varepsilon=(\det \boldsymbol{L}_\varepsilon)(\det \boldsymbol{A}_\varepsilon)$, we have   
\begin{equation}
	\label{eq:detAF}
 	(1/C)\det \boldsymbol{A}_\varepsilon \leq \det \boldsymbol{F}_\varepsilon \leq C\det \boldsymbol{A}_\varepsilon.
\end{equation} \EEE 
 Using the elementary estimate
\begin{equation*}
	(q+1)(h_q(t)+1)\geq q h_q(t) + (q+1)=t^q + \frac{ q}{t}\geq t^q \quad \text{for all $t> 0$}
\end{equation*}
together with \eqref{eq:E-bdd1}, 
we deduce
\begin{equation}
	\label{eq:detF-bdd}
	\int_\Omega (\det \boldsymbol{F}_\varepsilon)^q\,\d\boldsymbol{x}\leq C \int_\Omega (\det \boldsymbol{A}_\varepsilon)^q\,\d\boldsymbol{x}\leq C \left(   \int_\Omega  h_q(\det \boldsymbol{A}_\varepsilon)\,\d \boldsymbol{x} + 1 \right) \leq C.
\end{equation}
 Thus, $(\det \boldsymbol{F}_\varepsilon)_{\varepsilon>0}$ is bounded in $L^q(\Omega)$.
In view of \eqref{eq:y}, we also  know that    $\det \boldsymbol{F}_\varepsilon \to 1$ almost everywhere in $\Omega$ for a not relabeled subsequence. 
 Therefore,  \eqref{eq:det1} follows by   Vitali's convergence theorem and Urysohn's property.
\EEE 

\emph{Step 2 (Compactness of magnetizations).}
We fix an exponent $\alpha >0$ satisfying
\begin{equation}
	\label{eq:alpha}
	\frac{\beta q}{2(q-1)} < \alpha < 1,
\end{equation}
noting that such a choice is always possible in view of \eqref{eq:beta}.  Let  $ \Omega_\varepsilon\EEE \coloneqq \{ \boldsymbol{x}\in\Omega: |\boldsymbol{G}_\varepsilon(\boldsymbol{x})|<\varepsilon^{-\alpha} \}$. By Chebyshev's inequality, 
\begin{equation}
	\label{eq:chebyshev}
	\leb(\Omega \setminus   \Omega_\varepsilon\EEE)\leq \varepsilon^{2\alpha} \int_\Omega |D\boldsymbol{u}_\varepsilon|^2\,\d\boldsymbol{x}\leq C \varepsilon^{2\alpha} \quad \text{for all $\varepsilon >0$.}
\end{equation}  
Recall that
\begin{equation}
	\label{eq:detI}
	\det \left( \boldsymbol{I}+\varepsilon \boldsymbol{G} \right)=1+\sum_{k=1}^{3} \varepsilon^k P_k(\boldsymbol{G}) \quad \text{for all   $\boldsymbol{G}\in\Rtt$,}
\end{equation}
where each $P_k(\boldsymbol{G})$ is  a \EEE homogeneous polynomial of degree $k$ in the entries of $\boldsymbol{G}$ and, as such,
\begin{equation}
	\label{eq:detI2}
	|P_k(\boldsymbol{G})|\leq C |\boldsymbol{G}|^k \quad \text{for all  $k=1,2,3$ and $\boldsymbol{G}\in\Rtt$.}
\end{equation}
Thus,  recalling \eqref{eq:alpha}, we have
\begin{equation}
	\label{eq:det-A}
	\|\det \boldsymbol{F}_\varepsilon\|_{L^\infty(\Omega_\varepsilon\EEE)}\leq 1 + C \sum_{k=1}^3 \varepsilon^{k(1-\alpha)}\leq C \quad \text{for all $\varepsilon >0$.}
\end{equation}
By considering the Taylor expansion of $\Phi$ around $\boldsymbol{m}_\varepsilon$ as in \eqref{eq:taylor-Phi} and Proposition \ref{prop:img=imt}, we write
\begin{equation}
	\label{eq:lb-Phi}
	\begin{split}
	  \frac{1}{\varepsilon^\beta}\int_{\imt(\boldsymbol{y}_\varepsilon,\Omega)} \Phi & \left( \left( \boldsymbol{F}_\varepsilon^\top \circ \boldsymbol{y}_\varepsilon^{-1} \right)\boldsymbol{m}_\varepsilon  \right)\,\d\boldsymbol{\xi}\\
	 &\geq \frac{1}{\varepsilon^\beta} \int_{\img(\boldsymbol{y}_\varepsilon, \Omega_\varepsilon\EEE)} \Phi \left( \left( \boldsymbol{F}_\varepsilon^\top \circ \boldsymbol{y}_\varepsilon^{-1} \right)\boldsymbol{m}_\varepsilon  \right)\,\d\boldsymbol{\xi}\\
	 &= \frac{1}{\varepsilon^\beta} \int_{\img(\boldsymbol{y}_\varepsilon, \Omega_\varepsilon\EEE )} \Phi (\boldsymbol{m}_\varepsilon)\,\d\boldsymbol{\xi}\\
	  &+ \varepsilon^{1-\beta} \int_{\img(\boldsymbol{y}_\varepsilon, \Omega_\varepsilon\EEE)} D\Phi(\boldsymbol{m}_\varepsilon) \cdot \left(  \left ( \boldsymbol{G}_\varepsilon^\top \circ \boldsymbol{y}_\varepsilon^{-1} \right)\boldsymbol{m}_\varepsilon\right)\,\d\boldsymbol{\xi}\\
	  &+\frac{1}{\varepsilon^\beta} 	 	\int_{\img(\boldsymbol{y}_\varepsilon, \Omega_\varepsilon\EEE)} \eta_\Phi\left (\varepsilon \left ( \boldsymbol{G}_\varepsilon^\top \circ \boldsymbol{y}_\varepsilon^{-1} \right)\boldsymbol{m}_\varepsilon; \boldsymbol{m}_\varepsilon\right)\,\d\boldsymbol{\xi}.
	\end{split}
\end{equation} 
For the last two integrals on the right-hand side of \eqref{eq:lb-Phi}, we observe that 
\begin{equation*}
	\begin{split}
		\Bigg | \varepsilon^{1-\beta} \int_{\img(\boldsymbol{y}_\varepsilon, \Omega_\varepsilon\EEE)} D\Phi(\boldsymbol{m}_\varepsilon) \cdot &\left(  \left ( \boldsymbol{G}_\varepsilon^\top \circ \boldsymbol{y}_\varepsilon^{-1} \right)\boldsymbol{m}_\varepsilon\right)\,\d\boldsymbol{\xi} \Bigg | \\
		&\leq \varepsilon^{1-\beta} \|D\Phi\|_{L^\infty(\St;\Rt)} \int_{\img(\boldsymbol{y}_\varepsilon, \Omega_\varepsilon\EEE)}  \left |\boldsymbol{G}_\varepsilon^\top \circ \boldsymbol{y}_\varepsilon^{-1}\right | \,\d\boldsymbol{\xi}\\
		&\leq C  \varepsilon^{1-\beta} \int_{ \Omega_\varepsilon} |\boldsymbol{G}_\varepsilon|\,\det  \boldsymbol{F}_\varepsilon\,\d\boldsymbol{x}\leq C \varepsilon^{1-\beta}
	\end{split}
\end{equation*}
and
\begin{equation*}
	\begin{split}
		\Bigg |\frac{1}{\varepsilon^\beta} 	 	\int_{\img(\boldsymbol{y}_\varepsilon, \Omega_\varepsilon\EEE)} \eta_\Phi\Big (\varepsilon& \Big ( \boldsymbol{G}_\varepsilon^\top \circ \boldsymbol{y}_\varepsilon^{-1} \Big)\boldsymbol{m}_\varepsilon; \boldsymbol{m}_\varepsilon\Big)\,\d\boldsymbol{\xi} \Bigg | \\
		&\leq  \zeta_\Phi(\varepsilon^{1-\alpha}) \varepsilon^{1-\beta} \int_{\img(\boldsymbol{y}_\varepsilon, \Omega_\varepsilon\EEE)} |D\boldsymbol{G}_\varepsilon^\top\circ \boldsymbol{y}_\varepsilon^{-1}|\,\d\boldsymbol{\xi}\\
		&=  \zeta_\Phi(\varepsilon^{1-\alpha}) \varepsilon^{1-\beta} \int_{ \Omega_\varepsilon} |\boldsymbol{G}_\varepsilon| \det \boldsymbol{F}_\varepsilon\,\d\boldsymbol{x}\leq C  \varepsilon^{1-\beta}
	\end{split}
\end{equation*}
 thanks to the change-of-variable formula, \eqref{eq:zeta-Phi}, the boundedness of $(\boldsymbol{G}_\varepsilon)$ in $L^1(\Omega;\Rtt)$, and \eqref{eq:det-A}. Thus, from \eqref{eq:lb-Phi}, we deduce
\begin{equation}
	\label{eq:lb-Phi2}
	 \frac{1}{\varepsilon^\beta}\int_{\imt(\boldsymbol{y}_\varepsilon,\Omega)} \Phi  \left( \left( \boldsymbol{F}_\varepsilon^\top \circ \boldsymbol{y}_\varepsilon^{-1} \right)\boldsymbol{m}_\varepsilon  \right)\,\d\boldsymbol{\xi}\geq \frac{1}{\varepsilon^\beta} \int_{\img(\boldsymbol{y}_\varepsilon, \Omega_\varepsilon\EEE)} \Phi(\boldsymbol{m}_\varepsilon)\,\d\boldsymbol{\xi} - C\varepsilon^{1-\beta}. \EEE
\end{equation}
The integral on the right-hand side can be rewritten as
\begin{equation*}
	 \frac{1}{\varepsilon^\beta} \int_{\img(\boldsymbol{y}_\varepsilon, \Omega_\varepsilon\EEE)} \Phi (\boldsymbol{m}_\varepsilon)\,\d\boldsymbol{\xi}=\frac{1}{\varepsilon^\beta}\int_{\imt(\boldsymbol{y}_\varepsilon,\Omega)} \Phi(\boldsymbol{m}_\varepsilon)\,\d\boldsymbol{\xi} - \EEE\frac{1}{\varepsilon^\beta} \int_{\img(\boldsymbol{y}_\varepsilon,\Omega\setminus  \Omega_\varepsilon\EEE)} \Phi(\boldsymbol{m}_\varepsilon)\,\d\boldsymbol{\xi}.
\end{equation*}
Using  \eqref{eq:detF-bdd},  \eqref{eq:chebyshev},  the change-of-variable formula, and H\"{o}lder's inequality,  the second term can be estimated as
\begin{equation*}
	\begin{split}
		\left | \frac{1}{\varepsilon^\beta} \int_{\img(\boldsymbol{y}_\varepsilon,\Omega\setminus  \Omega_\varepsilon\EEE)} \Phi(\boldsymbol{m}_\varepsilon)\,\d\boldsymbol{\xi} \right | &\leq \frac{\|\Phi\|_{L^\infty(\St)}}{\varepsilon^\beta} \leb(\img(\boldsymbol{y}_\varepsilon,\Omega \setminus  \Omega_\varepsilon\EEE))=\frac{C}{\varepsilon^\beta} \int_{\Omega \setminus  \Omega_\varepsilon\EEE} \det  \boldsymbol{F}_\varepsilon\,\d\boldsymbol{x}\\
		&\leq \frac{C}{\varepsilon^\beta} \| \det  \boldsymbol{F}_\varepsilon\|_{L^q(\Omega)} \leb(\Omega \setminus  \Omega_\varepsilon\EEE)^{1/q'}\leq C \varepsilon^{\frac{2\alpha(q-1)}{q}-\beta}.
	\end{split}
\end{equation*}
 Thus, from \eqref{eq:lb-Phi2}, we obtain
\begin{equation}
	\label{eq:lb-Phi3}
	\frac{1}{\varepsilon^\beta}\int_{\imt(\boldsymbol{y}_\varepsilon,\Omega)} \Phi  \left( \left( \boldsymbol{F}_\varepsilon^\top \circ \boldsymbol{y}_\varepsilon^{-1} \right)\boldsymbol{m}_\varepsilon  \right)\,\d\boldsymbol{\xi}\geq \frac{1}{\varepsilon^\beta} \int_{\imt(\boldsymbol{y}_\varepsilon,\Omega)} \Phi(\boldsymbol{m}_\varepsilon)\,\d\boldsymbol{\xi} - C \varepsilon^{\gamma }, 
\end{equation}
 where $\gamma\coloneqq  \left (  \frac{2\alpha(q-1)}{q}\EEE\wedge 1 \right) -\beta$. Noting that $\gamma >0$ in view of \eqref{eq:beta} and \eqref{eq:alpha}, from
 \eqref{eq:E-bdd} and \eqref{eq:lb-Phi3}  we obtain 
\begin{equation}
	\label{eq:K}
	\sup_{\varepsilon >0 } \left\{ \frac{1}{\varepsilon^\beta} \int_{\imt(\boldsymbol{y}_\varepsilon,\Omega)} \Phi(\boldsymbol{m}_\varepsilon)\,\d \boldsymbol{\xi} + \varepsilon^\beta \int_{\imt(\boldsymbol{y}_\varepsilon,\Omega)} |D\boldsymbol{m}_\varepsilon|^2\,\d\boldsymbol{\xi} \right\}\leq C_0, 
\end{equation}
where $C_0>0$ is a constant depending on the supremum in \eqref{eq:E-bdd}.

Given \eqref{eq:y}, by applying \cite[Lemma 2.24 and Lemma 3.6(a)]{barchiesi.henao.moracorral}, we find a Lipschitz domain $U\subset \subset \Omega$ such that $U=\imt(\boldsymbol{id},U)\subset \subset \imt(\boldsymbol{y}_\varepsilon,\Omega)$ along a not relabeled subsequence. Hence, \eqref{eq:K} yields 
\begin{equation}
	\label{eq:bve}
	  \sup_{\varepsilon>0} \left \{\frac{1}{\varepsilon^\beta} \int_{U} \Phi(\boldsymbol{m}_\varepsilon)\,\d \boldsymbol{x} + \varepsilon^\beta\int_{U} |D\boldsymbol{m}_\varepsilon|^2\,\d\boldsymbol{x} \right \} \leq C_0.
\end{equation}
Let $i=1,\dots,M$ and, recalling \eqref{eq:fi}, define $w_{i,\varepsilon}\coloneqq f_i\circ \boldsymbol{m}_\varepsilon$.  By \cite[Lemma 4.2]{anzellotti.baldo.visintin}, we infer that $w_{i,\varepsilon} \in W^{1,2}(U)$ with
\begin{equation}
	\label{eq:D-bv}
	\text{$|D w_{i,\varepsilon}(\boldsymbol{x})|\leq \sqrt{\Phi(\boldsymbol{m}_\varepsilon(\boldsymbol{x}))}|D\boldsymbol{m}_\varepsilon(\boldsymbol{x})|$ \quad for almost all $\boldsymbol{x}\in U$.}	
\end{equation}
From \eqref{eq:bve}--\eqref{eq:D-bv} and Young's inequality, we find
\begin{equation*}
\sup_{\varepsilon >0 }\int_U |Dw_{i,\varepsilon}|\,\d\boldsymbol{x}\leq C_0,
\end{equation*}
while, trivially,
\begin{equation*}
	\sup_{\varepsilon>0}\|w_{i,\varepsilon}\|_{L^\infty(U)}\leq  \|f_i\|_{L^\infty(\St)}.
\end{equation*}
We deduce that the sequence $(w_{i,\varepsilon})_{\varepsilon>0}$ is bounded in $BV(U)$. Hence,  there exists $w_i\in BV(U)$ such that, up to subsequences,
\begin{equation}
	\label{eq:bv}
	w_{i,\varepsilon} \to w_i \quad \text{in $L^1(U)$,}
\end{equation}
as $\varepsilon \to 0^+$.
Since the constant $C_0$ in \eqref{eq:bve} does not depend on $U$, by applying \cite[Lemma 2.24 and Lemma 3.6(a)]{barchiesi.henao.moracorral} with a sequence of Lipschitz domains invading $\Omega$ together with a diagonal argument, we realize that $w_1,\dots,w_M\in BV(\Omega)$ and we select a not relabeled subsequence for which \eqref{eq:bv} holds true  for any domain $U\subset \subset \Omega$.

Let $E_i \coloneqq \{ \boldsymbol{x}\in\Omega:\: w_i(\boldsymbol{x})=0 \}$ and define $\boldsymbol{m}\colon \Omega \to \St$ by setting $\boldsymbol{m}\coloneqq \sum_{i=1}^M  \boldsymbol{b}_i \chi_{ E_i}$. We claim that 
\begin{equation}
	\label{eq:m-loc}
	\boldsymbol{m}_\varepsilon \to \boldsymbol{m} \quad \text{in $L^1(U)$} \quad \text{for any domain $U\subset \subset \Omega$.}
\end{equation}
Let $ F\coloneqq \{ \boldsymbol{x}\in U:\:\:\limsup_{\varepsilon \to 0^+} \Phi(\boldsymbol{m}_\varepsilon(\boldsymbol{x}))>0   \}$. From \eqref{eq:bve} and Fatou's lemma, we have 
\begin{equation*}
	\int_U \liminf_{\varepsilon \to 0^+} \Phi(\boldsymbol{m}_\varepsilon)\,\d \boldsymbol{x}\leq \liminf_{\varepsilon \to 0^+}  \int_U \Phi(\boldsymbol{m}_\varepsilon)\,\d \boldsymbol{x}=0
\end{equation*}
which shows that $\leb(F)=0$. Thus, for almost every $\boldsymbol{x}\in U$, the set of limit points of the sequence $(\boldsymbol{m}_\varepsilon(\boldsymbol{x}))_{\varepsilon>0}$ is contained in the set $ \{\boldsymbol{b}_1,\dots,\boldsymbol{b}_M  \}\EEE$. Let $i,j=1,\dots,M$ and recall \eqref{eq:fi}. For almost every $\boldsymbol{x}\in U\cap  E_i$, if $\boldsymbol{m}_\varepsilon(\boldsymbol{x})\to  \boldsymbol{b}_j\EEE$,  then $w_{i,\varepsilon}(\boldsymbol{x})=f_i(\boldsymbol{m}_\varepsilon(\boldsymbol{x}))\to f_i( \boldsymbol{b}_j\EEE)$ but also $w_{i,\varepsilon}(\boldsymbol{x})=f_i(\boldsymbol{m}_\varepsilon(\boldsymbol{x}))\to w_i(\boldsymbol{x})=0$,  as $\varepsilon\to 0^+$, \EEE in view of \eqref{eq:bv}. Thus, $f_i(  \boldsymbol{b}_j\EEE)=0$ which entails $i=j$. Therefore, $\boldsymbol{m}_\varepsilon(\boldsymbol{x})\to \boldsymbol{m}(\boldsymbol{x})$ for almost all $\boldsymbol{x}\in U$ and, by applying the dominated convergence theorem, we establish  \eqref{eq:m-loc}. Additionally, this argument  ensures that $ E_i\cap E_j=\emptyset$ for all $i\neq j$ as well as $\leb(\Omega \setminus \bigcup_{i=1}^M  E_i\EEE)=0$. Given the continuity of $f_i$, by applying once again  the dominated convergence theorem and comparing the result with \eqref{eq:bv}, we realize that $w_i=f_i\circ \boldsymbol{m}$ almost everywhere in $\Omega$. Therefore,  we deduce that each $ E_i$ has finite perimeter in $\Omega$ and, in turn,
 $\boldsymbol{m}\in BV(\Omega;\{  \boldsymbol{b}_1,\dots,\boldsymbol{b}_M  \})$. In particular, $(\boldsymbol{u},\boldsymbol{m})\in\mathcal{A}$. 
 Eventually, the combination of \eqref{eq:imt} and \eqref{eq:m-loc} easily gives \eqref{eq:m}. 

\emph{Step 3 (Convergence of compositions).} We move to the proof of \eqref{eq:composition}. From \eqref{eq:E-bdd1} and the trivial inequality
\begin{equation*}
	h_q(t)+\frac{q+1}{q}=\frac{t^q}{q}+\frac{1}{t}\geq \frac{1}{t} \quad \text{for all $t> 0$,}
\end{equation*}
we have
\begin{equation*}
	\int_\Omega \frac{1}{\det \boldsymbol{A}_\varepsilon}\,\d\boldsymbol{x}\leq  \int_\Omega h_q(\det \boldsymbol{A}_\varepsilon)\,\d\boldsymbol{x}+\frac{q+1}{q}\leb(\Omega)\leq C.
\end{equation*}
Then, thanks to   \eqref{eq:detAF},  we find 
\begin{equation*}
  	\int_\Omega \frac{1}{\det \boldsymbol{F}_\varepsilon}\,\d \boldsymbol{x}\leq C\int_\Omega \frac{1}{\det \boldsymbol{A}_\varepsilon}\,\d\boldsymbol{x}\leq C.
\end{equation*}
By  \cite[Remark 3.5]{bresciani.anisotropic}, we infer
  the equi-integrability of the sequence $(\chi_{\img(\boldsymbol{y}_\varepsilon,\Omega)}\det D \boldsymbol{y}^{-1}_\varepsilon)_{\varepsilon>0}$, while the same property trivially holds for $(\boldsymbol{m}_\varepsilon \circ \boldsymbol{y}_\varepsilon)_{\varepsilon>0}$. Therefore, given  \eqref{eq:y} and \eqref{eq:m}, the convergence in \eqref{eq:composition} is achieved  by applying   \cite[Proposition 3.4]{bresciani.friedrich.moracorral} and the dominated convergence theorem.
\end{proof}

\subsection{$\boldsymbol{\Gamma}$-convergence}

Next, we move to the proof of the $\Gamma$-convergence result. 
We will consider the second-order Taylor expansion of $W$ near the identity. Taking into account   (W2)--(W5),  we have
\begin{equation}
	\label{eq:taylor-W}
	W(\boldsymbol{I}+ \boldsymbol{B}\EEE)=\frac{1}{2} Q_W \EEE ( \boldsymbol{B})+ \eta_W( \boldsymbol{B}) \quad  \text{for all $\boldsymbol{B}\in \Rtt$,}
\end{equation}
where $ \eta_W (\boldsymbol{B})=\mathrm{o}(|\boldsymbol{B}|^2)$, as $|\boldsymbol{B}|\to 0^+$.   In particular, by (W3)--(W5), $W$ has a  \MMM global \EEE  minimum at the identity, so that $Q_W$ is a positive-semidefinite    quadratic form and, hence,  convex.   Additionally, in view of the minor symmetries of the elasticity tensor (see, e.g., \cite[p.~195]{gurtin}),
$\mathbf{C}_W(\boldsymbol{B})\in\mathrm{Sym}(3)$ and $\mathbf{C}_W(\boldsymbol{B})=\mathbf{C}_W(\sym \boldsymbol{B})$ for all $\boldsymbol{B}\in\Rtt$. Therefore, we have 
\begin{equation}
	\label{eq:Qsym}
	Q_W(\boldsymbol{B})=Q_W(\mathrm{sym}\boldsymbol{B}) \quad \text{for all $\boldsymbol{B}\in\Rtt$.}
\end{equation} 
 For $t>0$, \EEE  we define
\begin{equation*}
	\zeta_W( t \EEE)\coloneqq \sup \left\{  \frac{|\eta_W(\boldsymbol{B})|}{|\boldsymbol{B}|^2}:\:\:\boldsymbol{B}\in\Rtt,\:\:|\boldsymbol{B}|\leq  t \EEE   \right\},
\end{equation*}
so that assumption (W2) yields
\begin{equation}
	\label{eq:zeta-W}
	\lim_{ t \EEE\to 0^+} \zeta_W( t \EEE)=0.
\end{equation}

We begin by proving the lower bound.

\begin{proof}[Proof of claim (i) of Theorem \ref{thm:gamma}]
For convenience,  the proof is split into two steps.

\emph{Step 1 (Lower bound for the elastic energy).}
As in the proof of Theorem \ref{thm:compactness}, we set $	\boldsymbol{F}_\varepsilon\coloneqq D \boldsymbol{y}_\varepsilon$, $\boldsymbol{G}_\varepsilon \coloneqq D\boldsymbol{u}_\varepsilon$, and $\boldsymbol{L}_\varepsilon\coloneqq \boldsymbol{\Lambda}_\varepsilon(\boldsymbol{m}_\varepsilon \circ \boldsymbol{y}_\varepsilon)$. 
Additionally, we let  $\boldsymbol{K}_\varepsilon\coloneqq \boldsymbol{\Lambda}(\boldsymbol{m}_\varepsilon \circ \boldsymbol{y}_\varepsilon)$.
 Given \eqref{eq:Lambda-inverse},  we have
\begin{equation*}
	\boldsymbol{L}_\varepsilon^{-1}=\boldsymbol{I}-\varepsilon \boldsymbol{K}_\varepsilon+\varepsilon^2 \boldsymbol{J}_\varepsilon,
\end{equation*}
 where 
\begin{equation*}
	 	\boldsymbol{J}_\varepsilon\coloneqq \frac{a^2}{a_\varepsilon}  \boldsymbol{m}_\varepsilon \circ \boldsymbol{y}_\varepsilon \otimes \boldsymbol{m}_\varepsilon \circ \boldsymbol{y}_\varepsilon  + \frac{b^2}{b_\varepsilon} \left ( \boldsymbol{I}-\boldsymbol{m}_\varepsilon \circ \boldsymbol{y}_\varepsilon \otimes \boldsymbol{m}_\varepsilon \circ \boldsymbol{y}_\varepsilon\right ).
\end{equation*}
Similarly to \eqref{eq:ll}, we note that 
\begin{equation}
	\label{eq:JJ}
	 	\|\boldsymbol{J}_\varepsilon\|_{L^\infty(\Omega;\Rtt)}=\sqrt{\left( \frac{a^2}{a_\varepsilon}\right)^2+2\left (\frac{b^2}{b_\varepsilon}\right)^2  }\leq \sqrt{a^4+2b^4}\eqqcolon \tilde{\ell}.
\end{equation}
As \EEE  in Step 2 of the proof of Theorem \ref{thm:compactness}, define $\Omega_\varepsilon \EEE \coloneqq \{ \boldsymbol{x}\in\Omega:\:|\boldsymbol{G}_\varepsilon(\boldsymbol{x})|<\varepsilon^{-\alpha}  \}$ for  $ \alpha >0\EEE$   satisfying \EEE  \eqref{eq:alpha}. From \eqref{eq:u}, we recover \eqref{eq:chebyshev}. Thus,  $\chi_{A_\varepsilon}\to \Omega$ in measure and, from \eqref{eq:u} and \eqref{eq:composition}, we infer
\begin{align}
	\label{eq:u-A}
	\chi_{\Omega_\varepsilon\EEE}\boldsymbol{G}_\varepsilon \wk D\boldsymbol{u} \quad &\text{in $L^2(\Omega;\Rtt)$,}\\
	 \label{eq-K}
	 \text{$\chi_{\Omega_\varepsilon\EEE}\boldsymbol{K}_\varepsilon \to \boldsymbol{\Lambda}(\boldsymbol{m})$} \quad &\text{in $L^1(\Omega;\Rtt)$}.
\end{align} 
 For convenience, let $\boldsymbol{H}_\varepsilon \coloneqq \boldsymbol{J}_\varepsilon + (\varepsilon \boldsymbol{J}_\varepsilon - \boldsymbol{K}_\varepsilon) \boldsymbol{G}_\varepsilon$. From \eqref{eq:ll} and \eqref{eq:JJ}, we easily obtain
\begin{equation}
	\label{eq:HH}
	 \varepsilon \|\boldsymbol{H}_\varepsilon\|_{L^\infty(\Omega_\varepsilon\EEE;\Rtt)}\leq \varepsilon \tilde{\ell}+ \varepsilon^{1-\alpha} (\varepsilon \tilde{\ell}+\ell),
\end{equation}
where the right-hand side goes to zero  because of  \eqref{eq:alpha}.     
 Using the Taylor expansion \eqref{eq:taylor-W}, we write
\begin{equation*}
	\begin{split}
		\mathcal{E}_\varepsilon^{\rm e}(\boldsymbol{y}_\varepsilon,\boldsymbol{m}_\varepsilon)&\geq  \frac{1}{\varepsilon^2}\int_{ \Omega_\varepsilon\EEE} W(\boldsymbol{L}_\varepsilon^{-1}\boldsymbol{F}_\varepsilon)\,\d\boldsymbol{x}
		=\frac{1}{\varepsilon^2}\int_{ \Omega_\varepsilon \EEE} W\left( \boldsymbol{I}+\varepsilon (\boldsymbol{G}_\varepsilon-\boldsymbol{K}_\varepsilon)+\varepsilon^2 \boldsymbol{H}_\varepsilon  \right)\,\d\boldsymbol{x}\\
		&=\frac{1}{2} \int_{ \Omega_\varepsilon\EEE} Q_W(\boldsymbol{G}_\varepsilon -\boldsymbol{K}_\varepsilon +\varepsilon \boldsymbol{H}_\varepsilon)\,\d\boldsymbol{x} + \frac{1}{\varepsilon^2} \int_{ \Omega_\varepsilon } \eta_W \left( \varepsilon(\boldsymbol{G}_\varepsilon - \boldsymbol{K}_\varepsilon)+\varepsilon^2 \boldsymbol{H}_\varepsilon   \right)\,\d\boldsymbol{x}.
	\end{split}
\end{equation*} 
For the first integral on the right-hand side, owing to the convexity of $Q_W$  and \eqref{eq:Qsym}, \EEE we have 
\begin{equation*}
	\begin{split}
		\liminf_{\varepsilon \to 0^+} \int_{ \Omega_\varepsilon} Q_W\left( \boldsymbol{G}_\varepsilon - \boldsymbol{K}_\varepsilon+\varepsilon \boldsymbol{H}_\varepsilon \right)\,\d\boldsymbol{x}
		&\geq \int_\Omega Q_W(D\boldsymbol{u}-\boldsymbol{\Lambda}(\boldsymbol{m}))\,\d\boldsymbol{x}=\int_\Omega Q_W(E\boldsymbol{u}-\boldsymbol{\Lambda}(\boldsymbol{m}))\,\d\boldsymbol{x}
	\end{split}
\end{equation*}
thanks to \eqref{eq:u-A}--\eqref{eq:HH}. For the second one, we estimate
\begin{equation*}
	\begin{split}
		\bigg |\frac{1}{\varepsilon^2} \int_{ \Omega_\varepsilon\EEE } \eta_W&(\varepsilon (\boldsymbol{G}_\varepsilon-\boldsymbol{K}_\varepsilon)+\varepsilon^2 \boldsymbol{H}_\varepsilon)\,\d \boldsymbol{x} \bigg |\\
		&\leq \int_{ \Omega_\varepsilon\EEE} \zeta_W\left(\varepsilon \left ( |\boldsymbol{G}_\varepsilon|+ |\boldsymbol{K}_\varepsilon|\right )+\varepsilon^2 |\boldsymbol{H}_\varepsilon| \right) |\boldsymbol{G}_\varepsilon - \boldsymbol{K}_\varepsilon+\varepsilon \boldsymbol{H}_\varepsilon  |^2\,\d\boldsymbol{x} \\
		&\leq 3 \zeta_W\left (\varepsilon^{1-\alpha}+\varepsilon
		\ell +   \varepsilon \|\boldsymbol{H}_\varepsilon\|_{L^\infty( \Omega_\varepsilon;\Rtt)}\EEE\right ) \int_\Omega \left( |\boldsymbol{G}_\varepsilon|^2+|\boldsymbol{K}_\varepsilon|^2+|\varepsilon \boldsymbol{H}_\varepsilon|^2 \right)\,\d\boldsymbol{x}\\
		&\leq C \zeta_W\left (\varepsilon^{1-\alpha}+\varepsilon
		\ell +   \varepsilon \|\boldsymbol{H}_\varepsilon\|_{L^\infty( \Omega_\varepsilon;\Rtt)}\EEE \right),
	\end{split}
\end{equation*}
thanks to the boundedness of $(\boldsymbol{G}_\varepsilon)_{\varepsilon >0}$ in $L^2(\Omega;\Rtt)$ coming from \eqref{eq:u}, and to \eqref{eq:ll} and \eqref{eq:HH}. Hence, recalling \eqref{eq:alpha}, we obtain  
\begin{equation}
	\label{eq:lb-el}
	\liminf_{\varepsilon \to 0^+} \mathcal{E}^{\rm e}_\varepsilon(\boldsymbol{y}_\varepsilon,\boldsymbol{m}_\varepsilon)\geq \mathcal{E}^{\rm e}(\boldsymbol{u},\boldsymbol{m}).
\end{equation} \EEE 

\emph{Step 2 (Lower bound for the magnetic energy).} 
Without loss of generality, we assume \eqref{eq:E-bdd} and
we select a not relabeled subsequence for which the inferior limit of $(\mathcal{E}^{\rm m}_\varepsilon(\boldsymbol{y}_\varepsilon,\boldsymbol{m}_\varepsilon))_{\varepsilon>0}$ is attained as a limit. 
  
By arguing as in Step 2 of the proof of Theorem \ref{thm:compactness}, we recover \eqref{eq:lb-Phi3} for some $\gamma>0$ and, in turn,
\begin{equation}
	\label{eq:km}
	  \frac{1}{\varepsilon^\beta}\int_{\imt(\boldsymbol{y}_\varepsilon,\Omega)} \Phi(\boldsymbol{m}_\varepsilon)\,\d\boldsymbol{\xi}+\varepsilon^\beta \int_{\imt(\boldsymbol{y}_\varepsilon,\Omega)} |D\boldsymbol{m}_\varepsilon|^2\,\d\boldsymbol{\xi}  \leq  \mathcal{E}_\varepsilon^{\rm m}(\boldsymbol{y}_\varepsilon,\boldsymbol{m}_\varepsilon) + C \varepsilon^\gamma. 
\end{equation}  
Also, from \eqref{eq:E-bdd} and \eqref{eq:y}, we deduce that, up to subsequences,  $U=\imt(\boldsymbol{id},U)\subset \subset \imt(\boldsymbol{y}_\varepsilon,\Omega)$ for any domain $U\subset \subset \Omega$ by \cite[Lemma 2.24 and Lemma 3.6(a)]{barchiesi.henao.moracorral} and, hence, \eqref{eq:m-loc} holds true. Recalling \eqref{eq:fi}, again by \cite[Lemma 4.2]{baldo}, we have that $w_{i,\varepsilon}\coloneqq f_i\circ \boldsymbol{m}_\varepsilon\in W^{1,2}(U)$ satisfies \eqref{eq:D-bv} for all 
 $i=1,\dots,M$.  Using \eqref{eq:km} and  Young's inequality, we obtain \EEE 
\begin{equation*}
	 \int_E |Dw_{i,\varepsilon}|\,\d\boldsymbol{x}\leq
	\mathcal{E}_\varepsilon^{\rm m}(\boldsymbol{y}_\varepsilon,\boldsymbol{m}_\varepsilon) + C \varepsilon^\gamma
\end{equation*}
for any measurable subset $E\subset U$, and, in turn,
\begin{equation*}
	\left( \bigvee_{i=1}^M \left (|Dw_{i,\varepsilon}| \,\d\boldsymbol{x}  \EEE \right)  \right)(U)  = \EEE  \int_U  \max \{|Dw_{i,\varepsilon}|:\:i=1,\dots,M\}\EEE \,\d\boldsymbol{x} \leq 
	\mathcal{E}_\varepsilon^{\rm m}(\boldsymbol{y}_\varepsilon,\boldsymbol{m}_\varepsilon)+ C \varepsilon^\gamma,
\end{equation*}
where the equality is justified by \eqref{eq:vee}. \EEE 
From \eqref{eq:m-loc},   we obtain   \eqref{eq:bv} for $w_i\coloneqq f_i\circ \boldsymbol{m}$.  Therefore, taking the inferior limit, as $\varepsilon\to 0^+$, in the previous inequality by means of Lemma \ref{lem:tot-var}, then the supremum among all domains $U\subset \subset \Omega $,  and applying \cite[Proposition 2.2]{baldo},   we obtain
\begin{equation}
	\label{eq:lb-mag}
\liminf_{\varepsilon \to 0^+} 	\mathcal{E}_\varepsilon^{\rm m}(\boldsymbol{y}_\varepsilon,\boldsymbol{m}_\varepsilon)\geq \left( \bigvee_{i=1}^M |Dw_{i}| \right)(\Omega)=\mathcal{E}^{\rm m}(\boldsymbol{m}).
\end{equation}
Eventually, the combination of \eqref{eq:lb-el} and \eqref{eq:lb-mag} proves \eqref{eq:lb}.
\end{proof}

\begin{proof}[Proof of claim (ii) of  Theorem \ref{thm:gamma}]
We subdivide the proof into three steps.

\emph{Step 1 (Reduction to Lipschitz displacements).}
We argue that it is sufficient to prove the claim for Lipschitz displacements.  Given our assumptions on $\Delta$ and  $\boldsymbol{d}$, in view of \cite[Proposition A.2]{agostiniani.dalmaso.desimone}   there exists a sequence $(\boldsymbol{u}_n)_{n\in \N}$ in $W^{1,\infty}(\Omega;\Rt)$ such that
\begin{equation*}
	\boldsymbol{u}_n \to \boldsymbol{u} \quad \text{in $W^{1,2}(\Omega;\Rt)$,}
\end{equation*}
so that, by dominated convergence,
\begin{equation*}
	\lim_{n\to \infty} \mathcal{E}(\boldsymbol{u}_n,\boldsymbol{m})=\mathcal{E}(\boldsymbol{u},\boldsymbol{m}).
\end{equation*}
Hence, if we show that there exist a recovery sequence for each $(\boldsymbol{u}_n,\boldsymbol{m})$, then the claim for $\boldsymbol{u}$ can be established by means of a standard diagonal argument (see, e.g.,  \cite[Remark 1.29]{braides}). 
Henceforth, to avoid cumbersome  notation, we will simply assume that $\boldsymbol{u}\in W^{1,\infty}(\Omega;\Rt)$ and we will always consider its Lipschitz representative. 

\emph{Step 2 (Construction of the recovery sequence).}	
Setting $L \coloneqq \|D\boldsymbol{u}\|_{L^\infty(\Omega;\Rtt)}\vee 1$, we consider  $0<\varepsilon\leq ( 2L)^{-1}$. We define the Lipschitz map $\boldsymbol{y}_\varepsilon\coloneqq \boldsymbol{id}+\varepsilon \boldsymbol{u}$  and we check that it  is injective. Let $\boldsymbol{x}_1,\boldsymbol{x}_2\in\closure{\Omega}$ with $\boldsymbol{x}_1\neq \boldsymbol{x}_2$ and, by contradiction, suppose that $\boldsymbol{y}_\varepsilon(\boldsymbol{x}_1)=\boldsymbol{y}_\varepsilon(\boldsymbol{x}_2)$. Then
\begin{equation*}
	\begin{split}
		|\boldsymbol{x}_1-\boldsymbol{x}_2|&=|\boldsymbol{y}_\varepsilon(\boldsymbol{x}_1)-\boldsymbol{y}_\varepsilon(\boldsymbol{x}_2)-(\boldsymbol{x}_1-\boldsymbol{x}_2)|=\varepsilon |\boldsymbol{u}(\boldsymbol{x}_1)-\boldsymbol{u}(\boldsymbol{x}_2)|\\
		&\leq \varepsilon \,L|\boldsymbol{x}_1-\boldsymbol{x}_2|\leq \frac{1}{2} |\boldsymbol{x}_1-\boldsymbol{x}_2|
	\end{split}
\end{equation*} 
which provides a contradiction.  A similar argument shows that the inverse deformation $\boldsymbol{y}_\varepsilon^{-1}\colon \boldsymbol{y}_\varepsilon(\closure{\Omega}) \to \Rt$ is also Lipschitz: given $\boldsymbol{x}_1,\boldsymbol{x}_2\in\closure{\Omega}$, we have 
\begin{equation*}
	\begin{split}
		|\boldsymbol{x}_1-\boldsymbol{x}_2|&\leq |\boldsymbol{y}_\varepsilon(\boldsymbol{x}_1)-\boldsymbol{y}_\varepsilon(\boldsymbol{x}_2)|+\varepsilon |\boldsymbol{u}(\boldsymbol{x}_1)-\boldsymbol{u}(\boldsymbol{x}_2)|\\
		&\leq |\boldsymbol{y}_\varepsilon(\boldsymbol{x}_1)-\boldsymbol{y}_\varepsilon(\boldsymbol{x}_2)|+\varepsilon  L |\boldsymbol{x}_1-\boldsymbol{x}_2|
	\end{split}
\end{equation*}
so that 
\begin{equation*}
 |\boldsymbol{x}_1-\boldsymbol{x}_2|\leq \frac{1}{1-\varepsilon L} |\boldsymbol{y}_\varepsilon(\boldsymbol{x}_1)-\boldsymbol{y}_\varepsilon(\boldsymbol{x}_2)|\leq 2 |\boldsymbol{y}_\varepsilon(\boldsymbol{x}_1)-\boldsymbol{y}_\varepsilon(\boldsymbol{x}_2)|
\end{equation*}
since $\varepsilon \leq (2L)^{-1}$.
This proves the claim and shows that
\begin{equation*}
	\|D\boldsymbol{y}_\varepsilon^{-1}\|_{L^\infty(\boldsymbol{y}_\varepsilon(\Omega);\Rtt)}\leq  2 \quad \text{for all $\varepsilon \leq (2L)^{-1}$.}
\end{equation*}
 We observe that $\boldsymbol{y}\in\mathcal{Y}_{p,q}(\Omega)$,   $\imt(\boldsymbol{y}_\varepsilon,\Omega)=\boldsymbol{y}_\varepsilon(\Omega)$, and $\leb(\img(\boldsymbol{y},\Omega))=\leb(\boldsymbol{y}_\varepsilon(\Omega))$. \EEE  Also, \eqref{eq:y}--\eqref{eq:u} trivially hold.

 Let $\boldsymbol{m}\in BV(\Omega;\{ \boldsymbol{b}_1,\dots,\boldsymbol{b}_M  \})$.
 From \cite[Theorem 2.6]{anzellotti.baldo.visintin}, we infer the existence of a sequence $(\boldsymbol{\mu}_\varepsilon)_{\varepsilon>0}$ in $W^{1,2}(\Omega;\St)$ satisfying
\begin{equation}
	\label{eq:mu}
	\boldsymbol{\mu}_\varepsilon \to \boldsymbol{m} \quad \text{in $L^1(\Omega;\Rt)$,}
\end{equation}
and
\begin{equation}
	\label{eq:abv}
	\lim_{\varepsilon \to 0^+} \left\{ \frac{1}{\varepsilon^\beta} \int_\Omega \Phi(\boldsymbol{\mu}_\varepsilon)\,\d\boldsymbol{x}+\varepsilon^\beta \int_\Omega |D\boldsymbol{\mu}_\varepsilon|^2\,\d\boldsymbol{x}     \right\}=\mathcal{E}^{\rm m}(\boldsymbol{m}).
\end{equation}
We set $\boldsymbol{m}_\varepsilon\coloneqq \boldsymbol{\mu}_\varepsilon \circ \boldsymbol{y}_\varepsilon^{-1} \in W^{1,2}(\boldsymbol{y}_\varepsilon(\Omega);\St)$. To prove \eqref{eq:m}, fix $r>1$. 
By the area formula, 
\begin{equation*}
	\begin{split}
		\lim_{\varepsilon \to 0^+} \int_{\Rt}|\chi_{\boldsymbol{y}_\varepsilon(\Omega)}\boldsymbol{m}_\varepsilon|^r\,\d\boldsymbol{\xi}&= \lim_{\varepsilon \to 0^+}\leb(\boldsymbol{y}_\varepsilon(\Omega))= \lim_{\varepsilon \to 0^+}\int_\Omega \det D \boldsymbol{y}_\varepsilon \,\d\boldsymbol{x}\\
		&=\leb(\Omega)=\int_{\Rt} |\chi_\Omega \boldsymbol{m}|^r\,\d\boldsymbol{x},
	\end{split}
\end{equation*}
while, for any $\boldsymbol{\psi} \in C^0_{\rm c}(\Rt;\Rt)$, we have
\begin{equation*}
	\lim_{\varepsilon \to 0^+} \int_{\boldsymbol{y}_\varepsilon(\Omega)}\boldsymbol{m}_\varepsilon \cdot \boldsymbol{\psi}\,\d\boldsymbol{x}=\lim_{\varepsilon \to 0^+} \int_\Omega \boldsymbol{\mu}_\varepsilon \cdot \boldsymbol{\psi} \circ \boldsymbol{y}_\varepsilon\,\det D \boldsymbol{y}_\varepsilon \,\d\boldsymbol{x}=\int_\Omega \boldsymbol{m}\circ \boldsymbol{\psi}\,\d\boldsymbol{x}, 
\end{equation*}
thanks to the change-of-variable formula and the dominated convergence theorem with \eqref{eq:y}  and \eqref{eq:mu}.  Hence,  $\chi_{\boldsymbol{y}_\varepsilon(\Omega)}\boldsymbol{m}_\varepsilon \to \chi_\Omega \boldsymbol{m}$ in $L^r(\Rt;\Rt)$. As $\boldsymbol{y}_\varepsilon \to \boldsymbol{id}$ uniformly in $\Omega$, the set $\boldsymbol{y}_\varepsilon(\Omega)$ is contained in any compact neighborhood of $\Omega$, at least for $\varepsilon \ll 1$, so that the previous convergence implies  \eqref{eq:m}.  

\emph{Step 3 (Attainment of the lower bound).}
First, we prove that
\begin{equation}
	\label{eq:recover-el}
	\lim_{\varepsilon \to 0^+} \mathcal{E}_\varepsilon^{\rm e}(\boldsymbol{y}_\varepsilon,\boldsymbol{m}_\varepsilon)=\mathcal{E}^{\rm e}(\boldsymbol{u},\boldsymbol{m}).
\end{equation}
For convenience, noting that $\boldsymbol{m}_\varepsilon \circ \boldsymbol{y}_\varepsilon=\boldsymbol{\mu}_\varepsilon$, we set $\boldsymbol{F}_\varepsilon\coloneqq D\boldsymbol{y}_\varepsilon$, $\boldsymbol{L}_\varepsilon\coloneqq \boldsymbol{\Lambda}_\varepsilon(\boldsymbol{\mu}_\varepsilon)$, and $\boldsymbol{K}_\varepsilon\coloneqq \boldsymbol{\Lambda}(\boldsymbol{\mu}_\varepsilon)$.
Clearly, \eqref{eq:mu} yields
\begin{equation}
	\label{eq:M-recovery}
	\boldsymbol{K}_\varepsilon \to \boldsymbol{\Lambda}(\boldsymbol{m})\quad \text{in $L^2(\Omega;\Rtt)$}
\end{equation}
 by the dominated convergence theorem. Similarly to Step 1 in the proof of Theorem~\ref{thm:gamma}(i), we define
\begin{equation*}
 \boldsymbol{J}_\varepsilon\coloneqq \frac{a^2}{a_\varepsilon} \boldsymbol{\mu}_\varepsilon \otimes \boldsymbol{\mu}_\varepsilon +\frac{b^2}{b_\varepsilon}  \left(\boldsymbol{I}-\boldsymbol{\mu}_\varepsilon \otimes \boldsymbol{\mu}_\varepsilon \right)
\end{equation*}
which still fulfills \eqref{eq:JJ}. Also, we let $\boldsymbol{H}_\varepsilon\coloneqq \boldsymbol{J}_\varepsilon + (\varepsilon \boldsymbol{J}_\varepsilon - \boldsymbol{K}_\varepsilon)D\boldsymbol{u}$ and, similarly to \eqref{eq:HH},  we check that
\begin{equation}
	\label{eq:HH2}
	 \|\boldsymbol{H}_\varepsilon\|_{L^\infty(\Omega;\Rtt)}\leq \tilde{\ell}+(\varepsilon \tilde{\ell}+\ell)L  \leq \tilde{\ell}+(\tilde{\ell}+\ell)L\eqqcolon \hat{\ell}.
\end{equation}
By considering the Taylor expansion \eqref{eq:taylor-W}, we write
\begin{equation*}
	\begin{split}
		\mathcal{E}^{\rm e}_\varepsilon(\boldsymbol{y}_\varepsilon,\boldsymbol{m}_\varepsilon)&=\frac{1}{\varepsilon^2}\int_\Omega W(\boldsymbol{L}_\varepsilon^{-1}\boldsymbol{F}_\varepsilon)\,\d\boldsymbol{x}=\frac{1}{\varepsilon^2}\int_\Omega W(\boldsymbol{I}+\varepsilon(D\boldsymbol{u}-\boldsymbol{K}_\varepsilon)+\varepsilon^2 \boldsymbol{H}_\varepsilon)\,\d\boldsymbol{x}\\
		&= \frac{1}{2}\int_\Omega Q_W(D\boldsymbol{u}-\boldsymbol{K}_\varepsilon + \varepsilon \boldsymbol{H}_\varepsilon)\,\d \boldsymbol{x}+\frac{1}{\varepsilon^2}\int_\Omega \eta_W(\varepsilon(D\boldsymbol{u}-\boldsymbol{K}_\varepsilon)+\varepsilon^2 \boldsymbol{H}_\varepsilon)\,\d\boldsymbol{x}.
	\end{split}
\end{equation*}
Using  \eqref{eq:Qsym}, \EEE \eqref{eq:M-recovery}--\eqref{eq:HH2}, and the dominated convergence theorem, we compute
\begin{equation*}
	\lim_{\varepsilon \to 0^+}\frac{1}{2} \int_\Omega Q_W(D\boldsymbol{u}-\boldsymbol{K}_\varepsilon +\varepsilon \boldsymbol{H}_\varepsilon)\,\d\boldsymbol{x}=\frac{1}{2} \int_\Omega Q_W(D\boldsymbol{u}-\boldsymbol{\Lambda}(\boldsymbol{m}))\,\d\boldsymbol{x}=\mathcal{E}^{\rm e}(\boldsymbol{u},\boldsymbol{m}).
\end{equation*}
As
\begin{equation*}
	\begin{split}
		 \frac{1}{\varepsilon^2} \int_\Omega \eta_W \big( \varepsilon(D\boldsymbol{u}-\boldsymbol{K}_\varepsilon)+\varepsilon^2 \boldsymbol{H}_\varepsilon \big)\,\d\boldsymbol{x}\,\d\boldsymbol{x}&\leq \zeta_W(\varepsilon(L+\ell)+\varepsilon^2 \hat{\ell}) \int_\Omega |D\boldsymbol{u}-\boldsymbol{K}_\varepsilon+\varepsilon \boldsymbol{H}_\varepsilon|^2\,\d\boldsymbol{x}\\
		&\leq 3 \zeta_W(\varepsilon(L+\ell)+\varepsilon^2 \hat{\ell}\,) \int_\Omega \left( |D\boldsymbol{u}|^2+|\boldsymbol{K}_\varepsilon|^2+\varepsilon^2|\boldsymbol{H}_\varepsilon|^2  \right)\,\d\boldsymbol{x}\\
		&\leq 3  \left( L^2+\ell^2+\varepsilon^2 \hat{\ell}^{\,2}  \right) \leb(\Omega) \zeta_W(\varepsilon(L+\ell)+\varepsilon^2 \hat{\ell}\,)
	\end{split}
\end{equation*}
in view of \eqref{eq:ll} and \eqref{eq:HH2},
then \eqref{eq:recover-el} follows thanks to \eqref{eq:zeta-W}.

Next, we prove that
\begin{equation}
	\label{eq:recover-mag}
	\lim_{\varepsilon \to 0^+} \mathcal{E}^{\rm m}_\varepsilon(\boldsymbol{y}_\varepsilon,\boldsymbol{m}_\varepsilon)=\mathcal{E}^{\rm m}(\boldsymbol{m}).
\end{equation}
Having \eqref{eq:abv} in mind,  it is sufficient to check that 
\begin{equation}
	\label{eq:mag1}
	\lim_{\varepsilon \to 0^+}\frac{1}{\varepsilon^\beta} \left |  \int_{\boldsymbol{y}_\varepsilon(\Omega)} \Phi \left( (D\boldsymbol{y}_\varepsilon^\top \circ \boldsymbol{y}_\varepsilon^{-1}) \boldsymbol{m}_\varepsilon \right)\,\d\boldsymbol{\xi} - \int_\Omega \Phi(\boldsymbol{\mu}_\varepsilon)\,\d \boldsymbol{x}   \right |=0
\end{equation}
and
\begin{equation}
	\label{eq:mag2}
	\lim_{\varepsilon \to 0^+} \varepsilon^\beta \left | \int_{\boldsymbol{y}_\varepsilon(\Omega)} |D\boldsymbol{m}_\varepsilon|^2\,\d \boldsymbol{\xi} - \int_\Omega |D\boldsymbol{\mu}_\varepsilon|^2\,\d\boldsymbol{x} \right |=0.
\end{equation}

For  \eqref{eq:mag1}, recalling \eqref{eq:detI}--\eqref{eq:detI2}, we  observe that $\det \boldsymbol{F}_\varepsilon=1+\varepsilon g_\varepsilon$ for a  sequence $(g_\varepsilon)_{\varepsilon>0}$ bounded in $L^\infty(\Omega)$.  By applying the change-of-variable formula and considering the Taylor expansion of $\Phi$  in \eqref{eq:taylor-Phi} around $\boldsymbol{\mu}_\varepsilon$, we write
\begin{equation*}
	\begin{split}
			\int_{\boldsymbol{y}_\varepsilon(\Omega)}  &\Phi \left( (D\boldsymbol{y}_\varepsilon^\top \circ \boldsymbol{y}_\varepsilon^{-1}) \boldsymbol{m}_\varepsilon \right)\,\d\boldsymbol{\xi} - \int_\Omega \Phi(\boldsymbol{\mu}_\varepsilon)\,\d \boldsymbol{x}=\int_\Omega \left \{ \Phi \left( \boldsymbol{F}_\varepsilon^\top \boldsymbol{\mu}_\varepsilon \right)\det \boldsymbol{F}_\varepsilon - \Phi(\boldsymbol{\mu}_\varepsilon)\right \}\,\d\boldsymbol{x}\\
			&=\varepsilon \int_\Omega D\Phi(\boldsymbol{\mu}_\varepsilon)\cdot\left ( (D\boldsymbol{u})^\top \boldsymbol{\mu}_\varepsilon \right)\,\d\boldsymbol{x} + \int_\Omega \eta_\Phi \left( \varepsilon (D\boldsymbol{u})^\top \boldsymbol{\mu}_\varepsilon  \right)\,\d\boldsymbol{x}
			+\varepsilon \int_\Omega g_\varepsilon \Phi\left( \boldsymbol{F}_\varepsilon^\top \boldsymbol{\mu}_\varepsilon  \right)\,\d\boldsymbol{x}.
	\end{split}
\end{equation*} 
Thus, given a neighborhood $V$ of $\St$, for $\varepsilon \ll 1$ we find 
\begin{equation*}
	\begin{split}
		\frac{1}{\varepsilon^\beta} \Bigg |  \int_{\boldsymbol{y}_\varepsilon(\Omega)} \Phi & \Big( (D\boldsymbol{y}_\varepsilon^\top \circ \boldsymbol{y}_\varepsilon^{-1}) \boldsymbol{m}_\varepsilon \Big)\,\d\boldsymbol{\xi} - \int_\Omega \Phi(\boldsymbol{\mu}_\varepsilon)\,\d \boldsymbol{x}   \Bigg |\\
		&\leq \left \{ \varepsilon^{1-\beta} \left (L  \|D\Phi\|_{L^\infty(\St;\Rt)}+
		 \|g_\varepsilon\|_{L^\infty(\Omega)}\|\Phi\|_{L^\infty(V)}    \right)
		 +\varepsilon^{ 1-\beta\EEE} L\zeta_\Phi\left( \varepsilon L \right)\right\} \leb(\Omega),
	\end{split}
\end{equation*}
so that \eqref{eq:mag1} follows from \eqref{eq:beta} and \eqref{eq:zeta-Phi}.

For \eqref{eq:mag2}, we observe that $\boldsymbol{F}_\varepsilon^{-1}=\boldsymbol{I}+ \varepsilon \EEE\boldsymbol{E}_\varepsilon$ for a sequence $(\boldsymbol{E}_\varepsilon)_{\varepsilon>0}$ bounded in $L^\infty(\Omega;\Rtt)$.
Indeed, using Neumann series, 
\begin{equation*}
	\begin{split}
		\boldsymbol{F}_\varepsilon^{-1}&= \boldsymbol{I}+\sum_{k=1}^\infty (-\varepsilon)^k (D\boldsymbol{u})^k= \boldsymbol{I}+ \varepsilon \sum_{k=1}^\infty (-1)^k\varepsilon^{k-1}  (D\boldsymbol{u})^k=\boldsymbol{I}+\varepsilon \sum_{l=0}^\infty (-1)^{l+1}\varepsilon^l(D\boldsymbol{u})^{l+1},
	\end{split}
\end{equation*}
where,  for $\varepsilon\leq (2L)^{-1}$, 
\begin{equation*}
	\left | \sum_{l=0}^\infty (-1)^{l+1}\varepsilon^l  (D\boldsymbol{u})^{l+1} \right |\leq  L \sum_{l=0}^\infty |{\varepsilon}L|^l\leq L \sum_{l=0}^\infty \frac{1}{2^l}= 2L.
\end{equation*}
Additionally, $D\boldsymbol{m}_\varepsilon=(D\boldsymbol{\mu}_\varepsilon)\circ \boldsymbol{y}_\varepsilon^{-1}(D\boldsymbol{y}_\varepsilon)^{-1}\circ \boldsymbol{y}_\varepsilon^{-1}$ almost everywhere in $\boldsymbol{y}_\varepsilon(\Omega)$ thanks to the chain rule, and \eqref{eq:abv} entails
\begin{equation}
	\label{eq:recovery-mu}
	\sup_{\varepsilon >0 } \varepsilon^\beta  \int_\Omega |\mathbf{M}_\varepsilon|^2\,\d\boldsymbol{x}\leq C_0,
\end{equation}
where we set $\boldsymbol{\rm M}_\varepsilon\coloneqq D\boldsymbol{\mu}_\varepsilon$. By applying 
 the change-of-variable formula, we write
\begin{equation*}
	\begin{split}
		\int_{\boldsymbol{y}_\varepsilon(\Omega)} |D\boldsymbol{m}_\varepsilon|^2\,\d\boldsymbol{\xi}-\int_\Omega |D\boldsymbol{\mu}_\varepsilon|^2\,\d\boldsymbol{x}&=\int_\Omega \left \{ |\boldsymbol{\rm M}_\varepsilon\boldsymbol{F}_\varepsilon^{-1}|^2\det \boldsymbol{F}_\varepsilon - |\boldsymbol{\rm M}_\varepsilon|^2\right \}\,\d\boldsymbol{x}\\
		&= \varepsilon \int_\Omega \left( g_\varepsilon |\mathbf{M}_\varepsilon|^2+2 \left( \mathbf{M}_\varepsilon : (\mathbf{M}_\varepsilon \boldsymbol{E}_\varepsilon) \right) (1+\varepsilon g_\varepsilon)   \right) \,\d\boldsymbol{x}\\
		&+\varepsilon^2 \int_\Omega |\mathbf{M}_\varepsilon \boldsymbol{E}_\varepsilon|^2 (1+\varepsilon g_\varepsilon)\,\d\boldsymbol{x}
	\end{split}
\end{equation*}
Therefore, \eqref{eq:recovery-mu} and the boundedness of the sequences $(g_\varepsilon)_{\varepsilon >0}$ and $(\boldsymbol{E}_\varepsilon)_{\varepsilon >0}$ in $L^\infty(\Omega)$ and $L^\infty(\Omega;\Rtt)$, respectively, yield
\begin{equation*}
	\begin{split}
		\varepsilon^\beta \left | \int_{\boldsymbol{y}_\varepsilon(\Omega)} |D\boldsymbol{m}_\varepsilon|^2\,\d \boldsymbol{\xi} - \int_\Omega |D\boldsymbol{\mu}_\varepsilon|^2\,\d\boldsymbol{x} \right |&\leq C_0\left( \|g_\varepsilon\|_{L^\infty(\Omega)} + 2 \|\boldsymbol{E}_\varepsilon\|_{L^\infty(\Omega;\Rtt)}\,\|1+\varepsilon g_\varepsilon\|_{L^\infty(\Omega)}   \right) \varepsilon\\
		&+C_0 \|1+\varepsilon g_\varepsilon\|_{L^\infty(\Omega)}\,\|\boldsymbol{E}_\varepsilon\|^2_{L^\infty(\Omega;\Rtt)} \varepsilon^2,
	\end{split}
\end{equation*}
and, in turn, \eqref{eq:mag2} follows.

This concludes the proof of \eqref{eq:recover-mag} which, together with \eqref{eq:recover-el}, gives \eqref{eq:olb}.
\end{proof}

\subsection{Convergence of almost minimizers}

 Before presenting the proof of Corollary \ref{cor:am}, we discuss the well posedness of the steady Maxwell system \eqref{eq:maxwell}. More generally, given $s\in (1,\infty)$ and $\boldsymbol{\zeta}\in L^s(\Rt;\Rt)$, we consider the system
\begin{equation}
	\label{eq:maxwell-g}
	\begin{cases}
		\curl\,\boldsymbol{h}=\boldsymbol{0} \\
		\div\, \boldsymbol{h}=\div\, \boldsymbol{\zeta}
	\end{cases} \quad \text{in $\Rt$.}
\end{equation} 
Solutions  to  the previous system will be considered in the distributional sense.
Below, we refer to the homogeneous Sobolev space (sometimes named after Beppo Levi or Deny-Lions \cite{deny.lions}) defined as
\begin{equation*}
	V^{1,s}(\Rt)\coloneqq \left\{ v  \in W^{1,s}_{\rm loc}(\Rt):\:D  v  \in L^s(\Rt;\Rt)  \right\}.
\end{equation*}

\begin{proposition}[Well-posedness of the Maxwell system]
	\label{prop:maxwell}
Let $\boldsymbol{\zeta}\in L^s(\Rt;\Rt)$ with $s\in (1,\infty)$. Then,  there exists a unique   distributional solution $\boldsymbol{h}_{\boldsymbol{\zeta}}\in L^s(\Rt;\Rt)$ to \eqref{eq:maxwell-g}. Moreover, the solution operator $\boldsymbol{\zeta}\mapsto \boldsymbol{h}_{\boldsymbol{\zeta}}$  is linear and bounded on $L^s(\Rt;\Rt)$. 	
\end{proposition}
\begin{proof}
In view of the first equation, we can write $\boldsymbol{h}=-D v $ for some $ v \in V^{1,s}(\Rt)$. Hence,  the system \eqref{eq:maxwell-g} is equivalent to the single equation
\begin{equation*}
	\div \left( -D v  + \boldsymbol{\zeta}\right)=0 \quad \text{in $\Rt$.}
\end{equation*} 
By \cite[Theorem 5.1]{praetorius}, the latter admits a unique distributional solution $ v_{\boldsymbol{\zeta}}\EEE\in V^{1,s}(\Rt)$ up to additive constants. The same result shows that the map $\boldsymbol{\zeta}\mapsto \boldsymbol{h}_{\boldsymbol{\zeta}}\coloneqq D v_{\boldsymbol{\zeta}}\EEE$ defines a bounded linear operator on $L^s(\Rt;\Rt)$. 
\end{proof}

For the proof of Corollary \ref{cor:am}, we need the following result.

\begin{proposition}[Additional compactness]
	\label{prop:sf}
Under the same assumptions of Theorem \ref{thm:compactness}, we additionally have  
\begin{align}
	\label{eq:inv-det}
	\frac{1}{\det D \boldsymbol{y}_\varepsilon} \to 1 \quad &\text{in $L^1(\Omega)$,}\\
	\label{eq:mm}
	\chi_{\imt(\boldsymbol{y}_\varepsilon,\Omega)}\boldsymbol{m}_\varepsilon \det D \boldsymbol{y}_\varepsilon^{-1}\to  \chi_\Omega \boldsymbol{m} \quad &\text{in $L^2(\Rt;\Rt)$.}
\end{align}  
\end{proposition}
\begin{proof}
First, we prove \eqref{eq:inv-det}. For convenience, let again  $\boldsymbol{F}_\varepsilon\coloneqq D \boldsymbol{y}_\varepsilon$  and $\boldsymbol{A}_\varepsilon\coloneqq \boldsymbol{\Lambda}_\varepsilon^{-1}(\boldsymbol{m}_\varepsilon \circ \boldsymbol{y}_\varepsilon)\boldsymbol{F}_\varepsilon$. \EEE By (W5), we have 
\begin{equation*}
 \int_\Omega h_q( \det \boldsymbol{A}_\varepsilon \EEE)\,\d\boldsymbol{x}\leq \frac{1}{c_W} \int_\Omega    W(\boldsymbol{A}_\varepsilon)\,\d\boldsymbol{x} \EEE  \leq C\varepsilon^2,
\end{equation*}
so that $h_q(\det \boldsymbol{A}_\varepsilon)\to 0$ in $L^1(\Omega)$. Recalling that $\det \boldsymbol{F}_\varepsilon=c_\varepsilon \det \boldsymbol{A}_\varepsilon$ by \EEE  \eqref{eq:det-Lambda}, we estimate  
\begin{equation*}
	\begin{split}
	\int_\Omega h_q(\det \boldsymbol{F}_\varepsilon)\,\d\boldsymbol{x}&=\int_\Omega \left \{ \frac{c_\varepsilon^q (\det \boldsymbol{A}_\varepsilon)^q}{q}+\frac{1}{c_\varepsilon \det \boldsymbol{A}_\varepsilon} - \frac{q+1}{q}   \right \}\,\d\boldsymbol{x}\\
		&\leq \left( c_\varepsilon^q \vee c_\varepsilon^{-1} \right) \int_\Omega  \left \{ \frac{(\det \boldsymbol{A}_\varepsilon)^q}{q}+\frac{1}{ \det \boldsymbol{A}_\varepsilon} \right \}\,\d\boldsymbol{x} - \frac{q+1}{q}\leb(\Omega)\\ 
		&= \left( c_\varepsilon^q \vee c_\varepsilon^{-1} \right) \int_\Omega h_q(\det \boldsymbol{A}_\varepsilon)\,\d\boldsymbol{x}+ \left( \left(c_\varepsilon^q \vee c_\varepsilon^{-1} \right) -1  \right) \frac{q+1}{q}\leb(\Omega),
	\end{split}
\end{equation*}
\EEE where the right-hand side goes to zero as $c_\varepsilon \to 1$. Thus,
\begin{equation}
	\label{eq:hh}
	\text{ $h_q(\det \boldsymbol{F}_\varepsilon)\to 0$ in $L^1(\Omega)$.}
\end{equation}
As in Step 1 of Theorem \ref{thm:compactness}, we have \eqref{eq:det1} and,  from   the convergence in measure, we obtain
\begin{equation}
	\label{eq:conv-measure}
	\text{$\chi_{\{ \det \boldsymbol{F}_\varepsilon>\frac{q}{q+1} \}}\to 1$ in $L^1(\Omega)$}, \qquad \text{$\chi_{\{ \det \boldsymbol{F}_\varepsilon\leq \frac{q}{q+1} \}}\to 0$ in $L^1(\Omega)$.}
\end{equation}
 Using \MMM \eqref{eq:det1} and  dominated convergence \EEE    we easily see that
\begin{equation*}
	\lim_{\varepsilon \to 0^+} \int_{ \{ \det \boldsymbol{F}_\varepsilon>\frac{q}{q+1} \}} \left | \frac{1}{\det \boldsymbol{F}_\varepsilon}-1 \right |\,\d\boldsymbol{x} =0.
\end{equation*}
Then,
\begin{equation*}
	\begin{split}
		\int_{ \{ \det \boldsymbol{F}_\varepsilon\leq \frac{q}{q+1} \} } \left | \frac{1}{\det \boldsymbol{F}_\varepsilon}-1 \right |\,\d\boldsymbol{x}&=\int_{ \{ \det \boldsymbol{F}_\varepsilon\leq \frac{q}{q+1} \} } \left ( \frac{1}{\det \boldsymbol{F}_\varepsilon}-1 \right )\,\d\boldsymbol{x}\\
		&=\int_{ \{ \det \boldsymbol{F}_\varepsilon\leq \frac{q}{q+1} \} } \left ( \frac{1}{\det \boldsymbol{F}_\varepsilon}-\frac{q+1}{q} \right )\,\d\boldsymbol{x}+\frac{1}{q}\leb(\{ \det \boldsymbol{F}_\varepsilon\leq \textstyle  \frac{q}{q+1} \})\\
		&\leq \int_\Omega h_q(\det \boldsymbol{F}_\varepsilon)\,\d\boldsymbol{x}+\frac{1}{q}\leb(\{ \det \boldsymbol{F}_\varepsilon\leq \textstyle  \frac{q}{q+1} \}),
	\end{split}
\end{equation*}
where the right-hand side goes to zero in view of \eqref{eq:hh} and \eqref{eq:conv-measure}. Therefore, \eqref{eq:inv-det} follows. \EEE 

Next, we look at \eqref{eq:mm}. For convenience, we write 
 $\boldsymbol{\zeta}_\varepsilon\coloneqq \chi_{\imt(\boldsymbol{y}_\varepsilon,\Omega)}\boldsymbol{m}_\varepsilon\det D \boldsymbol{y}_\varepsilon^{-1}$.  By applying the change-of-variable formula, we compute
\begin{equation*}
	\int_{\Rt} |\boldsymbol{\zeta}_\varepsilon|^2\,\d\boldsymbol{\xi}=\int_{\imt(\boldsymbol{y}_\varepsilon,\Omega)} |\det D \boldsymbol{y}_\varepsilon^{-1}|^2\,\d\boldsymbol{\xi}=\int_\Omega \frac{1}{\det D \boldsymbol{y}_\varepsilon}\,\d \boldsymbol{x}.
\end{equation*}
Given \eqref{eq:inv-det}, letting $\varepsilon\to 0^+$ in the previous equation yields
\begin{equation}
	\label{eq:l2norm}
	\lim_{\varepsilon \to 0^+} \|\boldsymbol{\zeta}_\varepsilon\|^2_{L^2(\Rt;\Rt)}=\lim_{\varepsilon \to 0^+} \left \|  \frac{1}{\det D \boldsymbol{y}_\varepsilon}\right \|_{L^1(\Omega)}=\leb(\Omega)=\|\chi_\Omega \boldsymbol{m} \|^2_{L^2(\Rt;\Rt)}. 
\end{equation}
Let $\boldsymbol{\psi}\in C^0_{\rm c}(\Rt;\Rt)$. By employing once again the change-of-variable formula, with the aid of the dominated convergence theorem, using \eqref{eq:y} and \eqref{eq:m}--\eqref{eq:composition}, we obtain
\begin{equation*}
	\lim_{\varepsilon \to 0^+} \int_{\Rt} \boldsymbol{\zeta}_\varepsilon \cdot \boldsymbol{\psi}\,\d\boldsymbol{\xi}=\lim_{\varepsilon \to 0^+}\int_\Omega \boldsymbol{m}_\varepsilon \circ \boldsymbol{y}_\varepsilon \cdot \boldsymbol{\psi}\circ \boldsymbol{y}_\varepsilon\,\d\boldsymbol{x}=\int_\Omega \boldsymbol{m}\cdot \boldsymbol{\psi}\,\d\boldsymbol{x}.
\end{equation*}
Thus, $\boldsymbol{\zeta}_\varepsilon \wk \chi_\Omega \boldsymbol{m}$ in $L^2(\Rt;\Rt)$ which, together with \eqref{eq:l2norm}, gives \eqref{eq:mm}.
\end{proof}

\EEE 

With Proposition \ref{prop:sf} at hand,  the result for almost minimizers follows  from Theorem \ref{thm:compactness} and Theorem \ref{thm:gamma} in a standard manner.

\begin{proof}[Proof of Corollay \ref{cor:am}]
We consider a sequence   $(({\boldsymbol{w}}_\varepsilon,{\boldsymbol{n}}_\varepsilon))_{\varepsilon>0}$ with $({\boldsymbol{w}}_\varepsilon,{\boldsymbol{n}}_\varepsilon)\in\mathcal{A}_\varepsilon$ for all $\varepsilon>0$ such that 
\begin{equation*}
	\sup_{\varepsilon>0} \mathcal{E}({\boldsymbol{w}}_\varepsilon,{\boldsymbol{n}}_\varepsilon)<+\infty.
\end{equation*}
 By applying Theorem \ref{thm:gamma}(ii) for limit $(\boldsymbol{d},\boldsymbol{n})\in\mathcal{A}$ with    ${\boldsymbol{n}}\colon \Omega \to \St$  constantly equal to  $\boldsymbol{b}_1$,   we see that such sequences exist. 
  From \EEE Proposition \ref{prop:sf} to    $(({\boldsymbol{w}}_\varepsilon,{\boldsymbol{n}}_\varepsilon))_{\varepsilon>0}$, \EEE we infer that the sequence $({\boldsymbol{\zeta}}_\varepsilon)_{\varepsilon>0}$ with  ${\boldsymbol{\zeta}}_\varepsilon\coloneqq \chi_{\imt({\boldsymbol{w}}_\varepsilon,\Omega)}{\boldsymbol{n}}_\varepsilon \det D {\boldsymbol{w}}_\varepsilon^{-1}$ \EEE
 is bounded in $L^2(\Rt;\Rt)$. Thus, 
 \begin{equation*}
 	\sup_{\varepsilon>0} \mathcal{H}( {\boldsymbol{w}}_\varepsilon,{\boldsymbol{n}}_\varepsilon\EEE )<+\infty
 \end{equation*}
in view of Proposition \ref{prop:maxwell}. Using H\"{o}lder's inequality, we immediately find
\begin{equation*}
	\sup_{\varepsilon>0}|\mathcal{F}( {\boldsymbol{w}}_\varepsilon,{\boldsymbol{n}}_\varepsilon\EEE)|\leq \|\boldsymbol{f}\|_{L^1(\Rt;\Rt)}.	
\end{equation*}
Hence,
\begin{equation*}
	\sup_{\varepsilon>0} G_\varepsilon\leq \sup_{\varepsilon>0}  \mathcal{G}_\varepsilon( {\boldsymbol{w}}_\varepsilon,{\boldsymbol{n}}_\varepsilon\EEE)<+\infty,
\end{equation*}
 so that \eqref{eq:am} entails
 \begin{equation*}
 	\sup_{\varepsilon>0} \mathcal{G}_\varepsilon(\boldsymbol{y}_\varepsilon,\boldsymbol{m}_\varepsilon)<+\infty.
 \end{equation*}
  From this, we easily get
 \begin{equation*}
 	\sup_{\varepsilon>0} \mathcal{E}(\boldsymbol{y}_\varepsilon,\boldsymbol{m}_\varepsilon)\leq \sup_{\varepsilon>0} \left \{ \mathcal{G}_\varepsilon(\boldsymbol{y}_\varepsilon,\boldsymbol{m}_\varepsilon)+\mathcal{F}(\boldsymbol{y}_\varepsilon,\boldsymbol{m}_\varepsilon) \right \}<+\infty.
 \end{equation*}
By applying Theorem \ref{thm:compactness} and Proposition \ref{prop:sf} to the sequence $((\boldsymbol{y}_\varepsilon,\boldsymbol{m}_\varepsilon))_{\varepsilon >0}$, we identify $(\boldsymbol{u},\boldsymbol{m})\in \mathcal{A}$ for which \eqref{eq:y}--\eqref{eq:composition} and  \eqref{eq:inv-det}--\eqref{eq:mm} hold true.  Proposition \ref{prop:maxwell} ensures that $ \boldsymbol{h}_\varepsilon\coloneqq \EEE \boldsymbol{h}^{\boldsymbol{y}_\varepsilon,\boldsymbol{m}_\varepsilon}$ and $ \boldsymbol{h} \coloneqq \EEE \boldsymbol{h}^{\boldsymbol{id},\boldsymbol{m}}$ are well-defined and belong to $L^2(\Rt;\Rt)$. In particular, linearity and boundedness of the solution operator of the Maxwell system yield the estimate
\begin{equation*}
	\| \boldsymbol{h}_\varepsilon - \boldsymbol{h} \EEE \|_{L^2(\Rt;\Rt)}\leq C \| \chi_{\imt(\boldsymbol{y},\Omega)}\boldsymbol{m}_\varepsilon\det D \boldsymbol{y}_\varepsilon^{-1}-\chi_\Omega \boldsymbol{m} \|_{L^2(\Rt;\Rt)} 
\end{equation*}
and, in turn,
\begin{equation*}
	\label{eq:lb-H}
	\lim_{\varepsilon \to 0^+} \mathcal{H}(\boldsymbol{y}_\varepsilon,\boldsymbol{m}_\varepsilon)=\mathcal{H}(\boldsymbol{id},\boldsymbol{m})
\end{equation*}
thanks to \eqref{eq:mm}. 
Having in mind \eqref{eq:m}, we obtain
\begin{equation*}
	\label{eq:lb-F}
	\lim_{\varepsilon\to 0^+}\mathcal{F}(\boldsymbol{y}_\varepsilon,\boldsymbol{m}_\varepsilon)=\mathcal{F}(\boldsymbol{id},\boldsymbol{m}).
\end{equation*}
At this point, the minimality of $(\boldsymbol{u},\boldsymbol{m})$ and \eqref{eq:GGG} follow from Theorem \ref{thm:gamma} by standard \hbox{$\Gamma$-convergence} arguments (see, e.g., \cite[Remark 1.7]{braides}).
\end{proof}

\subsection{Regularity  of  the limiting problem}
 Recalling that for $M=2$ we have \eqref{eq:due}, our problem is related to the regularity of 
perimeter minimizers. \EEE  Specifically, we will refer to the notion of  sets   of almost minimal boundary   \cite{tamanini}. 

\begin{definition}[Set of almost minimal boundary]
Let $\varkappa\in (0,1)$. A set $E\in\mathcal{P}(\Rt)$ has $\varkappa$-almost minimal boundary in $\Omega$ if, for every set $A\subset \subset \Omega$, there exist $c_A,R_A>0$ such that, for every  ${\boldsymbol{a}}\in A$, \EEE $\rho\in \left( 0, R_A \right)$,  and   $F\in \mathcal{P}(\Rt)$ with $E\triangle {F} \subset \subset B({\boldsymbol{a}},\rho)$, \EEE there holds
\begin{equation*}
	 	\per (E;B({\boldsymbol{a}},\rho))\leq \per ( {F};B({\boldsymbol{a}},\rho))\EEE +c_A \rho^{2(1+\varkappa)}.
\end{equation*}
\end{definition}

In view of \cite[Theorem 1]{tamanini}, we known that $\partial^* E \cap \Omega=\partial E \cap \Omega$  is a surface of class $C^{1,\varkappa}$ whenever $E$ has $\varkappa$-almost minimimal boundary in $\Omega$.

We now present the proof of the regularity result. For convenience, we fix $\omega_3\coloneqq \frac{4\pi}{3}$.

\begin{proof}[Proof of Theorem \ref{thm:regularity}]
We subdivide the proof into three steps.	

\emph{Step 1 (Regularity of the displacement).}	
Let $ {\boldsymbol{v}}\EEE \in W^{1,2}(\Omega;\Rt)$ with $ {\boldsymbol{v}}\EEE=\boldsymbol{d}$ on $\Delta$. 
Testing the minimality of $( \boldsymbol{u},\boldsymbol{m}\EEE)$ by taking $( {\boldsymbol{v}},\boldsymbol{m} \EEE)\in\mathcal{A}$ as a competitor, we realize that $ \boldsymbol{u} \EEE$ is a minimizer of $ {\boldsymbol{v}}\mapsto \mathcal{G}({\boldsymbol{v}},\boldsymbol{m})\EEE$. Thus, $ \boldsymbol{u}=\boldsymbol{d}+\hat{\boldsymbol{u}}\EEE$, where $ \hat{\boldsymbol{u}}\EEE\in W^{1,2}(\Omega;\Rt)$ is the unique weak solution of the boundary value problem
\begin{equation*}
	\begin{cases}
		-\div \left( \mathbf{C}_W(E\hat{\boldsymbol{u}})  \right)=-\div \left( \mathbf{C}_W \left( \boldsymbol{\Lambda}( \boldsymbol{m}\EEE)-E\boldsymbol{d}   \right) \right) & \text{in $\Omega$,}\\
		 \hat{\boldsymbol{u}}\EEE=\boldsymbol{0} & \text{on $\Delta$,}\\
		 \mathbf{C}_W(E \hat{\boldsymbol{u}})\EEE\boldsymbol{\nu}_\Omega=-\mathbf{C}_W(E\boldsymbol{d})\boldsymbol{\nu}_\Omega & \text{on $\partial \Omega \setminus \Delta$.}
	\end{cases}
\end{equation*}
Here, $\mathbf{C}_W\colon \Rtt \to \Rtt$ and $\boldsymbol{\nu}_\Omega \colon \partial \Omega \to \St$   denote the elasticity tensor of $W$ and the outer unit normal to $\Omega$, respectively.
Note that $\mathbf{C}_W \left( \boldsymbol{\Lambda}( \boldsymbol{m}\EEE)-E \boldsymbol{d}   \right) \in L^\infty(\Omega;\Rtt)$. Therefore, standard elliptic regularity theory (see, e.g., \cite[Theorem 16.4]{ambrosio.pde}) yields $ \hat{\boldsymbol{u}}\EEE\in W^{1,r}(\Omega;\Rt)$ and, in turn, $ \boldsymbol{u}\EEE\in W^{1,r}(\Omega;\Rt)$ for all $r\in [1,\infty)$. 

\emph{Step 2 (Regularity of the jump set).}
Let $ {E} \EEE\coloneqq \{  {\boldsymbol{m}} =\boldsymbol{b}_1 \EEE  \}$ and $\varkappa\coloneqq {(s_{\rm{f}}-3)}/({2s_{\rm{f}}})\in (0,{1}/{2})$. We claim that $ E \EEE$ is a set of $\varkappa$-almost minimimal boundary in $\Omega$. To see this, we consider  $A\subset \subset \Omega$ and we set $R_A\coloneqq \dist(A;\partial \Omega)\wedge 1$. Let $ {\boldsymbol{a}} \EEE\in A$, $\rho\in (0,R_A)$, and $ {F} \EEE\in \mathcal{P}(\Rt)$ with $ E\triangle {F} \EEE \subset \subset B( {\boldsymbol{a}}\EEE,\rho)$. 
We need to show that
\begin{equation}
	\label{eq:amb}
	\per ( E\EEE;B( {\boldsymbol{a}}\EEE,\rho))\leq \per ( {F} \EEE;B( {\boldsymbol{a}}\EEE,\rho))+c_A\rho^{2(1+\varkappa)}
\end{equation}  
for some constant $c_A>0$ to be determined.
Define $ {\boldsymbol{n}}\EEE\in BV(\Omega;\{  \boldsymbol{b}_1,\boldsymbol{b}_2 \EEE \})$ by setting $ {\boldsymbol{n}}\EEE\coloneqq  \boldsymbol{b}_1\EEE$ in $ {F}$ and $ {\boldsymbol{n}}\EEE\coloneqq  \boldsymbol{b}_2 \EEE$ in $\Omega \setminus  {F}$. Then, $ {\boldsymbol{n}}=\boldsymbol{m}$ in $\Omega \setminus B( {\boldsymbol{a}},\rho \EEE)$. By minimality, $\mathcal{G}( \boldsymbol{u},\boldsymbol{m}\EEE)\leq \mathcal{G}( \boldsymbol{u}, {\boldsymbol{n}}\EEE)$ from which we obtain
\begin{equation}
	\label{eq:p}
	\begin{split}
		\per ( E;B({\boldsymbol{a}}\EEE,\rho))&\leq \per ( {F};B({\boldsymbol{a}}\EEE ,\rho))+\left (\mathcal{E}^{\rm e}( \boldsymbol{u},{\boldsymbol{n}})-\mathcal{E}^{\rm e}( \boldsymbol{u},\boldsymbol{m}\EEE) \right) \\
		&+\lambda \Big( \mathcal{H}(\boldsymbol{id}, {\boldsymbol{n}}\EEE)-\mathcal{H}(\boldsymbol{id}, \boldsymbol{m}\EEE) \Big)+\Big (\mathcal{F}(\boldsymbol{id}, {\boldsymbol{n}}\EEE)-\mathcal{F}(\boldsymbol{id}, \boldsymbol{m}\EEE) \Big) .
	\end{split}
\end{equation}
We estimate the three differences on the right-hand side of \eqref{eq:p}. We have
\begin{equation*}
	\begin{split}
		 \mathcal{E}^{\rm e}(\boldsymbol{u},{\boldsymbol{n}})-\mathcal{E}^{\rm e}(\boldsymbol{u},\boldsymbol{m})\EEE&=\frac{1}{2}\int_\Omega \left(  Q_W(\boldsymbol{\Lambda}( {\boldsymbol{n}}\EEE))-Q_W(\boldsymbol{\Lambda}( \boldsymbol{m}\EEE)) \right)\,\d\boldsymbol{x}\\
		&+\int_\Omega \mathbf{C}_W(E \boldsymbol{u}\EEE):\left(\boldsymbol{\Lambda}( \boldsymbol{m}\EEE)-\boldsymbol{\Lambda}( {\boldsymbol{n}}\EEE)  \right)\,\d\boldsymbol{x}.
	\end{split}
\end{equation*}
Recalling \eqref{eq:ll}, the first summand on the right-hand side is easily estimated as
\begin{equation}
	\label{eq:R1}
	\left | \int_\Omega \left(  Q_W(\boldsymbol{\Lambda}( {\boldsymbol{n}}\EEE))-Q_W(\boldsymbol{\Lambda}( \boldsymbol{m}\EEE)) \right)\,\d\boldsymbol{x} \right |\leq 2 \max \{ Q_W(\boldsymbol{E}):\:|\boldsymbol{E}|=\ell  \} \omega_3 \rho^3=C \rho^3,
\end{equation}
while for the second one we find 
\begin{equation}
	\label{eq:R2}
	\begin{split}
		\left | \int_\Omega \mathbf{C}_W(E \boldsymbol{u}\EEE):\left(\boldsymbol{\Lambda}( \boldsymbol{m}\EEE)-\boldsymbol{\Lambda}( {\boldsymbol{n}}\EEE)  \right)\,\d\boldsymbol{x}  \right |&\leq 2 \ell|\mathbf{C}_W|\int_{B( {\boldsymbol{a}}\EEE,\rho)} |E \boldsymbol{u}\EEE|\,\d\boldsymbol{x}\\
		&\leq 2 \ell  |\mathbf{C}_W| \|E \boldsymbol{u}\EEE \|_{L^r(\Omega;\Rtt)} \omega_3^{1/r'} \rho^{3/r'} \leq \EEE C \rho^{3/r'}
	\end{split}
\end{equation}
for any $r\in [1,\infty)$ thanks to H\"{o}lder's inequality. Writing $ {\boldsymbol{k}}\EEE\coloneqq \boldsymbol{h}^{\boldsymbol{id}, {\boldsymbol{n}}\EEE}$ and $ \boldsymbol{h}\EEE\coloneqq \boldsymbol{h}^{\boldsymbol{id}, \boldsymbol{m}\EEE}$, the second difference on  the right-hand side of \eqref{eq:p} is bounded by
\begin{equation*}
	\begin{split}
	 \left | 	\mathcal{H}(\boldsymbol{id}, {\boldsymbol{n}}\EEE)-\mathcal{H}(\boldsymbol{id}, \boldsymbol{m}\EEE)\right | &= \left | \int_{\Rt} | {\boldsymbol{k}}\EEE|^2\,\d\boldsymbol{\xi}-\int_{\Rt} | \boldsymbol{h}\EEE|^2\,\d\boldsymbol{\xi}\right |
		=\left |\int_{\Rt} ({\boldsymbol{k}}-\boldsymbol{h}\EEE)\cdot ( {\boldsymbol{k}}+\boldsymbol{h}\EEE)\,\d\boldsymbol{\xi} \right |\\
		&\leq \|   {\boldsymbol{k}}-\boldsymbol{h} \EEE \|_{L^s(\Rt;\Rt)}\,\left( \| {\boldsymbol{k}}\EEE\|_{L^{s'}(\Rt;\Rt)}+  \| \boldsymbol{h}\EEE\|_{L^{s'}(\Rt;\Rt)}  \right)
	\end{split}
\end{equation*}
for any $s \in (1,\infty)$. 
Then, from Proposition \ref{prop:maxwell}, we infer the estimates
\begin{equation*}
	\begin{split}
		\|{\boldsymbol{k}}-\boldsymbol{h}\EEE\|_{L^s(\Rt;\Rt)}&\leq C \|\chi_\Omega  {\boldsymbol{m}}\EEE-\chi_\Omega   \boldsymbol{m} \EEE\|_{L^s(\Rt;\Rt)}=C\| \chi_{B( {\boldsymbol{a}}\EEE,\rho)}({\boldsymbol{n}}-\boldsymbol{m}\EEE) \|_{L^s(\Rt;\Rt)}\\
		&\leq 2^{{1}/{s}}C\omega_3^{1/s}\rho^{{3}/{s}} \leq \EEE C\rho^{{3}/{s}} 
	\end{split}
\end{equation*}
and
\begin{equation*}
	\begin{split}
	\| {\boldsymbol{k}}\EEE\|_{L^{s'}(\Rt;\Rt)}+\| \boldsymbol{h} \EEE\|_{L^{s'}(\Rt;\Rt)}&\leq C \left( \| {\boldsymbol{n}}\EEE\|_{L^{s'}(\Omega;\Rt)} + \| \boldsymbol{m}\EEE\|_{L^{s'}(\Omega;\Rt)}   \right)
	=2C \leb(\Omega)^{1/s'} \leq \EEE C.	
	\end{split}
\end{equation*}
Thus,
\begin{equation}
	\label{eq:R3}
	\left | 	\mathcal{H}(\boldsymbol{id}, {\boldsymbol{n}}\EEE)-\mathcal{H}(\boldsymbol{id}, \boldsymbol{m}\EEE)\right | \leq C \rho^{3/s}.
\end{equation}
Eventually, for the last difference on the right-hand side of \eqref{eq:p}, we have
\begin{equation}
	\label{eq:R4}
	\begin{split}
		\left | \mathcal{F}(\boldsymbol{id}, {\boldsymbol{n}}\EEE) - \mathcal{F}(\boldsymbol{id}, \boldsymbol{m}\EEE)  \right |&\leq \int_\Omega |\boldsymbol{f}|\,| {\boldsymbol{n}}-\boldsymbol{m}\EEE|\,\d \boldsymbol{x}\leq 2 \int_{B( {\boldsymbol{a}}\EEE,\rho)} |\boldsymbol{f}|\,\d\boldsymbol{x}\\
		&\leq 2 \|\boldsymbol{f}\|_{L^{s_{\rm f}}(\Rt;\Rt)} \omega_3^{1/s_{\rm f}'} \rho^{3/s_{\rm f}'}\leq C \rho^{3/s_{\rm f}'}.
	\end{split}
\end{equation}
Altogether, \eqref{eq:p}--\eqref{eq:R4} yield the estimate 
\begin{equation*}
	\per ( E;B({\boldsymbol{a}}\EEE ,\rho))\leq \per ( {F};B({\boldsymbol{a}}\EEE,\rho))+ C \left( \rho^{3/r'}+\rho^{3/s}+\rho^{3/s_{\rm f}'} \right).
\end{equation*}
In particular, choosing   $r=s_{\rm f}$ and $s=s_{\rm f}'$, by noting that $3/s'_{\rm f}=2(1+\varkappa)$, we  obtain the desired estimate \eqref{eq:amb} with $c_A= C$. 
Therefore,   $ E\EEE$ is a set of $\varkappa$-almost minimimal boundary in $\Omega$, so that $J_{\boldsymbol{m}}=\partial^*E \cap \Omega$ is a surface of class $C^{1,\varkappa}$ in $\Omega$. \EEE

\emph{Step 3 (Improved regularity for the displacement).} Observe that $J_{ \boldsymbol{m}\EEE}$ is a closed subset of $\Omega$. Let $ U$ be one of the connected components of $\Omega \setminus J_{ \boldsymbol{m}\EEE}$. As $ \boldsymbol{m}\EEE$ is constant on $ U$, the displacement $ \hat{\boldsymbol{u}}=\boldsymbol{u}-\boldsymbol{d}\EEE$ is a weak solution to the equation
\begin{equation*}
	-\div \left( \mathbf{C}_W(E \hat{\boldsymbol{u}}\EEE) \right)=\div \left(\mathbf{C}_W(E\boldsymbol{d})\right) \quad \text{in $ U$.}
\end{equation*}
Hence, if $\boldsymbol{d}\in W^{2,s_{\rm d}}_{\rm loc}(\Omega;\Rt)$ with $s_{\rm d}>3$,  elliptic regularity theory (see, e.g., \cite[Remark 15.16]{ambrosio.pde}) yields $ \hat{\boldsymbol{u}}\EEE\in W^{2,s_{\rm d}}_{\rm loc}( U\EEE;\Rt)$ and, by Morrey's embedding, $ \boldsymbol{u}=\hat{\boldsymbol{u}}\EEE+\boldsymbol{d}\in C^{1,1-3/s_{\rm {d}}}( U \EEE;\Rt)$. As the connected component $ U \EEE$ is arbitrary, we conclude $ \boldsymbol{u}\EEE\in C^{1,1-3/s_{\rm {d}}}(\Omega \setminus J_{ \boldsymbol{m}\EEE};\Rt)$.
\end{proof}

\EEE 

\section*{Acknowledgements}

M.~B. acknowledges the support of the Alexander von Humboldt Foundation. He is member of GNAMPA (Gruppo Nazionale per l'Analisi Matematica, la Probabilit\'{a} e le loro Applicazioni) of INdAM (Istituto Nazionale di Alta Matematica).  The work of M.~F. was   supported
by the DAAD (Deutscher Akademischer Austauschdienst) project 57702972. \EEE

{
	}


\begin{thebibliography}{50}
	
\bibitem{agostiniani.dalmaso.desimone}
  {\sc V. Agostiniani, G. Dal Maso, A. De Simone}: \emph{Linear elasticity obtained from finite elasticity by $\Gamma$-convergence under weak coerciveness conditions}, Ann. Inst. H. Poincar\'{e} Anal. Non Lin\'{e}aire 29, 715--735 (2021).  
	
\bibitem{almi.kruzik.molchanova}
{\sc S. Almi, M. Kru\v{z}\'{i}k, A. Molchanova}: \emph{Linearization in magnetoelasticity}, Adv. Calc. Var. {18}, no. 2, 577--591 (2025).	
	
\bibitem{ambrosio.pde}
{\sc L. Ambrosio, A. Carlotto, A. Massaccesi:} \emph{Lecture Notes on Elliptic Partial Differential Equations}, Scuola Normale Superiore di Pisa (2010). Available at: \url{https://hdl.handle.net/11384/92283}. 
	

\bibitem{ambrosio.fusco.pallara}
{\sc L. Ambrosio, N. Fusco, D. Pallara}: \emph{Functions of Bounded Variations and Free Discontinuity Problems}, Oxford Mathematical Monographs, Oxford University Press (2000).

\bibitem{anzellotti.baldo.visintin}
{\sc G. Anzellotti, S. Baldo, A. Visintin}: \emph{Asymptotic Behavior of the Landau-Lifschitz Model of Ferromagnetism}, Appl. Math. Optim. 23, 171--192 (1991).

\bibitem{actuators}
{\sc V. Apicella, C.~S. Clemente, D. Davino, D. Leone, C. Visone}: \emph{Review of Modeling and Control of Magnetostrictive Actuators}, Actuators 8, Article no. 45 (2019).

\bibitem{baldo}
{\sc S. Baldo}: \emph{Minimal interface criterion for phase transitions 	in mixtures of Cahn-Hilliard fluids}, Anal. Inst. Henri Poincar\'{e} Analyse Non Lin\'{e}aire 7, no. 2, 67--90 (1990).

\bibitem{ball.gi}
{\sc J.~M.~Ball}: \emph{Global invertibility of Sobolev functions and the interpenetration	of matter},
  Proc. Roy. Soc. Edinburgh {88A}, 315--328  (1981).


\bibitem{ball.1982}
{\sc J. M. Ball}: \emph{Discontinuous equilibrium solutions and cavitation in nonlinear elasticity}, Phil. Trans. R. Soc. London A {306}, 557--611  (1982). \EEE



\bibitem{barchiesi.desimone}
{\sc M. Barchiesi, A. De Simone}: 
\newblock{\em Frank energy for nematic elastomers: a nonlinear model},  
\newblock{ESAIM Control Optim. Calc.
	Var.} {21}, no. 2, 372--377  (2015). \EEE

\bibitem{barchiesi.henao.moracorral}
{\sc M. Barchiesi, D. Henao, C. Mora-Corral}: \emph{Local Invertibility in Sobolev Spaces with
	Applications to Nematic Elastomers and
	Magnetoelasticity}, Arch. Ration. Mech. Anal. 224, 743--816  (2017).
	

	
\bibitem{braides}
{\sc A. Braides}: \emph{$\Gamma$-convergence for beginners}, Oxford Lecture Series in Mathematics and its Applications 22, Oxford University Press (2002). 	
	
\bibitem{bresciani.anisotropic}
{\sc M. Bresciani}: \emph{Anisotropic energies for the  modeling of cavitation in nonlinear elasticity}, submitted (2025). Preprint available at \url{https://arxiv.org/pdf/2503.23143}.


\bibitem{bresciani}
{\sc M. Bresciani}: \emph{Quasistatic evolutions in magnetoelasticity under subcritical coercivity assumptions}, Calc. Var. PDE 62, Art. no. 181 (2023).	

\bibitem{bresciani.m3as}
{\sc M. Bresciani}: \emph{Linearized von K\'{a}rm\'{a}n theory for incompressible magnetoelastic plates}, Math. Mod. Meth. Appl. Sci. 31, no. 10, 1987--2037 (2021).




\bibitem{bresciani.davoli.kruzik}
{\sc M. Bresciani, E. Davoli, M. Kru\v{z}\'{i}k}: \emph{Existence results in large-strain magnetoelasticity}, Ann. Inst. Henri Poincar\'{e} Anal. Non Lin\'{e}aire 40, no. 3 (2023).

\bibitem{bresciani.friedrich.moracorral}
{\sc M. Bresciani, M. Friedrich, M. Mora-Corral}: \emph{Variational models with Eulerian-Lagrangian formulation allowing for material failure},  Arch. Ration. Mech. Anal. {249}, no. 4 (2025).

\bibitem{bk}
{\sc M. Bresciani, M. Kru\v{z}\'{i}k}: \emph{A reduced model for plates arising as low energy $\Gamma$-limit in nonlinear magnetoelasticity}, SIAM J. Appl. Math. 55, no. 4, 3108--3168 (2023). 

\bibitem{bresciani.stroffolini}
{\sc M. Bresciani, B. Stroffolini}: \emph{Quasistatic evolution of Orlicz-Sobolev nematic elastomers}, Ann. Mat. Pura Appl. (2025). Article available at: \url{https://doi.org/10.1007/s10231-025-01580-1}.

\bibitem{brown}
{\sc W.~F.~Brown}: \emph{Magnetoelastic interactions}, Springer Tracts in Natural Philosophy 9, Springer-Verlag (1966). 

%

\bibitem{dalmaso}
{\sc G. Dal Maso}: \emph{An introduction to $\Gamma$-convergence}, Progress in Nonlinear Differential Equations and their Applications 8, Brikh\"{a}user (1993).

\bibitem{dalmaso.negri.percivale}
{\sc G. Dal Maso, M. Negri, D. Percivale}: \emph{Linearized Elasticity as $\Gamma$-limit of Finite Elasticity}, Set-Valued Var. Anal. 10, 165–183 (2002).


\bibitem{deny.lions}
{\sc J. Deny, J.-L- Lions}: \emph{Les espaces du type de Beppo Levi}, Ann. Inst. Fourier 5, 305--370 (1954).

\bibitem{desimone.podioguidugli}
{\sc A.~De Simone, P.~Podio-Guidugli}: \emph{On the Continuum Theory of Deformable Ferromagnetic Solids}, Arch. Ration. Mech. Anal. 136, 201--233 (1996).


%

\bibitem{engdhal}
{\sc G. Engdahl}: \emph{Handbook of Giant Magnetostrictive Materials}, Academic Press (2000).





\bibitem{fonseca.gangbo}
{\sc I. Fonseca, W. Gangbo}: \emph{Degree theory in analysis and applications}, Oxford Lecture Series in Mathematics and its Applications 2, Oxford University Press (1995).


\bibitem{fjm}
{\sc G. Friesecke, R. D. James, S. M\"{u}ller}: \emph{
A theorem on geometric rigidity and the derivation
of nonlinear plate theory from three-dimensional elasticity}, Comm. Pure Appl.
Math. 55, 1461–1506 (2002).


%
%
%
%
%

\bibitem{gkms-arma}
{\sc D. Grandi, M. Kru\v{z}\'{i}k, E. Mainini, U. Stefanelli}: \emph{A Phase-Field Approach to Eulerian Interfacial Energies}, Arch. Ration. Mech. Anal. 234, 351--373 (2019).

\bibitem{gkms}
{\sc D. Grandi, M. Kru\v{z}\'{i}k, E. Mainini, U. Stefanelli}: \emph{Equilibirum of Multiphase Solids with Eulerian Interfaces}, J. Elast. 142, 409--431 (2020).



\bibitem{sensors}
{\sc C.~A. Grimes, S.~C. Roy, S. Rani, Q. Cai}: \emph{Theory, Instrumentation and Applications of Magnetoelastic Resonance Sensors: a Review}, Sensors 11, 2809--2844 (2011).

\bibitem{gurtin}
 {\sc M.~E.~Gurtin}: \emph{An Introduction to Continuum Mechanics}, Mathematics in Sciences and Engineering 158, Academic Press (1981). \EEE 

%

\bibitem{henao.moracorral.invertibility}
{\sc D. Henao, C. Mora-Corral}: \emph{Invertibility and Weak Continuity of the Determinant for the Modelling of Cavitation and Fracture in Nonlinear Elasticity}, Arch. Ration. Mech. Anal. {197} (2010), 619--655.

\bibitem{henao.moracorral.fracture}
{\sc D. Henao, C. Mora-Corral}: \emph{Fracture surfaces and the regularity of inverses for $BV$
	deformations}, {Arch. Ration. Mech. Anal.} {201}, no. 2, 575--629  (2011). \EEE


%
%
%

\bibitem{henao.stroffolini}{\sc D. Henao, B. Stroffolini}: \emph{On Sobolev-Orlicz nematic elastomers}, {Nonlinear Anal.} {194}, 111513 (2020). 

\bibitem{hubert.schaefer}
{\sc A. Guber, R. Sch\"{a}fer}: \emph{Magnetic domains. The analysis of magnetic microstructures}, Springer (1998).

\bibitem{james.desimone}
 {\sc R.~D.~James, A.~De Simone}: \emph{A constrained theory of magnetoelasticity}, J. Mech. Phys. Solids 50, no. 2, 283--320 (2002).  

\bibitem{james.kinderlehrer}
{\sc R.~D.~James, D. Kinderlehrer}: \emph{Theory of magnetostriction with applications to $\mathrm{Tb}_x\mathrm{Dy}_{1-x}\mathrm{Fe}_2$}, Philos. Mag. B. 68, 237--274 (1993).


\bibitem{kruzik.stefanelli.zeman}
{\sc M. Kru\v{z}\'{i}k, U. Stefanelli, J. Zeman}: \emph{Existence results for incompressible magnetoelasticity}, Discrete Contin. Dyn. Syst.
35, 5999--6013   (2015). \EEE


\bibitem{landau.lifshitz}
{\sc L.~D.~Landau, E.~M.~Lifschitz}: \emph{Electrodynamics of continuous media}, Volume 8 of Course of Theoretical Physics, Pergamon Press (1984).

\bibitem{moracorral.oliva}
{\sc C. Mora-Corral, M. Oliva}: \emph{Relaxation of nonlinear elastic energies involving the deformed configuration and applications to nematic elastomers}, ESAIM  Control Optim. Calc. Var. 25, Article no. 19 (2019).

%
%
%
%
%
%
%
%
%
%





\bibitem{mueller.wcd}
{\sc  S. M\"{u}ller}: \emph{Weak continuity of determinants and nonlinear elasticity}, C. R. Acad. Sci. Paris S\'{e}r. I Math. 307, no. 9, 501--506
(1988).

\bibitem{mueller.det}
{\sc  S. M\"{u}ller}: \emph{$\mathrm{Det} = \mathrm{det}$. A Remark on the Distributional Determinant}, C. R.
Acad. Sci. Paris  S\'{e}r. I Math. 311, no. 1, 13--17 (1990).



%
%


\bibitem{praetorius}
{\sc D. Praetorius}: \emph{Analysis of the Operator $\Delta^{-1}\div$ Arising in Magnetic Models}, Z. Anal. Anwend. 23, no. 3, 589--605  (2004). \EEE


%
%
%

%
 
 \bibitem{scilla.stroffolini}
  {\sc G. Scilla, B. Stroffolini}: \EEE \emph{Relaxation of nonlinear elastic energies related to Orlicz–Sobolev nematic elastomers}, Rend. Lincei Mat. Appl. 31, 349--389 (2020).


\bibitem{tamanini}
{\sc I. Tamanini}: \emph{Boundaries of Caccioppoli sets with H\"{o}lder-continuous normal vector}, J. Reine Angew. Math. 334, 27--39  (1982). \EEE

%


\end{thebibliography}
\end{document}